\documentclass[11pt]{article}
\usepackage{amsfonts,amsmath}
\usepackage{latexsym}
\usepackage{amsthm}
\usepackage{amscd}
\usepackage{commath}
\usepackage{epsfig}
\usepackage{amssymb}
\usepackage{tikz}
\usepackage{graphicx}
\usepackage{caption}
\usepackage{enumitem}
\usepackage{hyperref}
\usepackage{mathtools}
\usepackage[toc,page]{appendix}
\newcommand{\tree}{\begin{tikzpicture}
		\draw[line width=0.06cm, black] (0,0) -- (0,1)  {};
		\draw[line width=0.06cm, black] (0,1) -- (1,2)  {};
		\draw[line width=0.06cm, black] (0,1) -- (0,2)  {};
		\draw[line width=0.06cm, black] (0,1) -- (-1,2)  {};
		\draw[line width=0.06cm, black] (-1,2) -- (-1.8,3)  {};
		\draw[line width=0.06cm, black] (-1,2) -- (-0.6,3)  {};
		\draw[line width=0.06cm, black] (1,2) -- (1,3)  {};
		\draw[line width=0.06cm, black] (1,2) -- (1.7,3)  {};
		\draw[line width=0.05cm, purple] (-1,2) -- (-2.2,3)  {};
		\draw[line width=0.05cm, purple] (-2.2,3) -- (-3,4)  {};
		\draw[line width=0.05cm, purple] (-2.2,3) -- (-2.4,4)  {};
		\draw[line width=0.05cm, blue] (-1,2) -- (-1.4,3)  {};
		\draw[line width=0.05cm, blue] (-1,2) -- (-0.95,3)  {};
		\draw[line width=0.05cm, blue] (-0.95,3) -- (-1.3,4)  {};
		\draw[line width=0.05cm, blue] (-0.95,3) -- (-0.8,4)  {};
		\draw[line width=0.05cm, blue] (-1.3,4) -- (-1.9,5)  {};
		\draw[line width=0.05cm, blue] (-1.3,4) -- (-1.3,5)  {};
		\draw[line width=0.05cm, blue] (-1.3,4) -- (-0.8,5)  {};
		\draw[line width=0.05cm, brown] (-1,2) -- (-0.1,3)  {};
		\draw[line width=0.05cm, brown] (-1,2) -- (0.4,3)  {};
		\draw[line width=0.05cm, brown] (0.4,3) -- (0,4)  {};
		\draw[line width=0.05cm, brown] (0.4,3) -- (0.6,4)  {};
		\node [fill=black, inner sep=2.5 pt,label=182: $i$] at (-1,2) {} ;
		\node [fill=black, inner sep=2.5 pt,label=left: $\phi_2(i)$] at (0,1) {} ;
		\node [fill=black, inner sep=2.5 pt,label=182: $i_0$] at (0,0) {} ;
\end{tikzpicture}}

\newcommand{\treeleaf}{\begin{tikzpicture}
		\node [fill=black, inner sep=2.5 pt,label=182: $i_0$] at (0,0) {} ;
		\draw[line width=0.05cm, black] (0,0) -- (0,1)  {};
		\node [fill=black, inner sep=2.5 pt,label=182: $i_1$] at (-1.8,3) {} ;
		\node [fill=black, inner sep=2.5 pt,label=182: $i_2$] at (-0.6,3) {} ;
		\node [fill=black, inner sep=2.5 pt,label=182: $i_3$] at (1,3) {} ;
		\node [fill=black, inner sep=2.5 pt,label=0: $i_4$] at (1.7,3) {} ;
		\node [fill=black, inner sep=1.5 pt] at (0,2) {} ;
		\node [fill=black, inner sep=1.5 pt] at (-1,2) {} ;
		\node [fill=black, inner sep=1.5 pt] at (1,2) {} ;
		\draw[line width=0.05cm, black] (0,1) -- (1,2)  {};
		\draw[line width=0.05cm, black] (0,1) -- (0,2)  {};
		\draw[line width=0.05cm, black] (0,1) -- (-1,2)  {};
		\draw[line width=0.05cm, black] (-1,2) -- (-1.8,3)  {};
		\draw[line width=0.05cm, black] (-1,2) -- (-0.6,3)  {};
		\draw[line width=0.05cm, black] (1,2) -- (1,3)  {};
		\draw[line width=0.05cm, black] (1,2) -- (1.7,3)  {};
		\draw [line width=0.04 cm, red] plot [smooth cycle] coordinates {(-1.8,3) (-2.2,5) (-1.3,5)};
		\draw [line width=0.04 cm, red] plot [smooth cycle] coordinates {(-0.6,3) (-0.9,4.5)  (-0.3,4.5)};
		\draw [line width=0.04 cm, red] plot [smooth cycle] coordinates {(1,3) (0.4,5.5) (1.4,5.5)};
		\draw [line width=0.04 cm, red] plot [smooth cycle] coordinates {(1.7,3) (1.7,5) (2.4,5)};
		\node[label=$T_1$] at (-1.83,3.7) {};
		\node[label=$T_2$] at (-0.6,3.65) {};
		\node[label=$T_3$] at (0.95,4.2) {};
		\node[label=$T_4$] at (1.95,4) {};
\end{tikzpicture}}

\newcommand{\treespine}{\begin{tikzpicture}
		\draw[line width=0.07cm, black] (0,0) -- (0,1)  {};
		\draw[line width=0.06cm, black] (0,1) -- (1,2)  {};
		\draw[line width=0.06cm, black] (0,1) -- (-1,2)  {};
		\draw[line width=0.06cm, black] (-1,2) -- (-1.8,3)  {};
		\draw[line width=0.06cm, black] (-1,2) -- (-0.6,3)  {};
		\draw[line width=0.06cm, black] (1,2) -- (0.7,3)  {};
		\draw[line width=0.05cm, black] (1,2) -- (1.7,3)  {};
		\draw [line width=0.04 cm, red] plot [smooth cycle] coordinates {(-1.8,3) (-2,7) (-1.6,7)};
		\draw [line width=0.04 cm, red] plot [smooth cycle] coordinates {(-0.6,3) (-0.8,7)  (-0.4,7)};
		\draw [line width=0.04 cm, red] plot [smooth cycle] coordinates {(0.7,3) (0.5,7) (0.9,7)};
		\draw [line width=0.04 cm, red] plot [smooth cycle] coordinates {(1.7,3) (1.5,7) (1.9,7)};
		\draw [line width=0.04 cm, blue] plot [smooth cycle] coordinates {(-1,2) (-3.5,3.4) (-2.7,3.4) (-1.55,2.4)};
		\draw [line width=0.04 cm, blue] plot [smooth cycle] coordinates {(-1,2) (-1.4,4.5) (-0.9,4.5) (-1,3)};
		\draw [line width=0.04 cm, blue] plot [smooth cycle] coordinates {(-1,2) (-0.5,2.6) (-0.25,3.3) (-0.15,3.3) (-0.2,2.5)};
		\draw [line width=0.04 cm, blue] plot [smooth cycle] coordinates {(1,2) (0.5,2.6) (0.3,3.3) (0.1,3.3) (0.2,2.5)};
		\draw [line width=0.04 cm, blue] plot [smooth cycle] coordinates {(1,2) (1,2.5) (1.1,4) (1.4,4)};
		\draw [line width=0.04 cm, blue] plot [smooth cycle] coordinates {(1,2) (1.7,2.7) (3,5) (2.5,3)};
		\draw [line width=0.04 cm, blue] plot [smooth cycle] coordinates {(0,1) (-0.2,2.2) (0.3,2.1) (0.05,1.2) };
		\draw [line width=0.04 cm, blue] plot [smooth cycle] coordinates {(0,1) (-3.3,2.4) (-2.1,2.2) };
		\draw [line width=0.04 cm, blue] plot [smooth cycle] coordinates {(0,1) (1.7,1.9) (2.4,2.2) (1.6,1.5)};
		\node [fill=black, inner sep=2 pt] at (-1.8,3) {} ;
		\node [fill=black, inner sep=2 pt] at (-0.6,3) {} ;
		\node [fill=black, inner sep=2 pt] at (0.7,3) {} ;
		\node [fill=black, inner sep=2 pt] at (1.7,3) {} ;
		\node [fill=black, inner sep=2 pt] at (-1,2) {} ;
		\node [fill=black, inner sep=2 pt] at (1,2) {} ;
		\node [fill=black, inner sep=2 pt] at (0,1) {} ;
		\node [fill=black, inner sep=2.5 pt,label=182: $i_0$] at (0,0) {} ;
\end{tikzpicture}}

\usepackage{float}

\usepackage{color,soul}

\newcommand{\bb}{\begin{equation}}
	\newcommand{\ee}{\end{equation}}

\newcommand{\half}{\frac{1}{2}}

\newcommand{\cT}{\mathcal T}
\newcommand{\cB}{\mathcal B}
\newcommand{\cF}{\mathcal F}

\newcommand{\D}{\mathbb{D}}

\definecolor{dblue}{rgb}{.61,.61,1}
\definecolor{lightblue}{rgb}{.61,.61,1}
\definecolor{altblue}{rgb}{.61,.61,1}

\tikzset{cross/.style={cross out, draw=black, minimum size=2*(#1-\pgflinewidth), inner sep=0pt, outer sep=0pt},
	cross/.default={1pt}}

\newtheorem{thm}{Theorem}[section]
\newtheorem{prop}[thm]{Proposition}
\newtheorem{lem}[thm]{Lemma}
\newtheorem{cor}[thm]{Corollary}
\newtheorem*{thm*}{Theorem}
\theoremstyle{remark}
\newtheorem{remark}{Remark}[section]

\theoremstyle{definition}

\def\IB{\relax\hbox{$\inbar\kern-.3em{\rm B}$}}
\def\IC{\relax\hbox{$\inbar\kern-.3em{\rm C}$}}
\def\ID{\relax\hbox{$\inbar\kern-.3em{\rm D}$}}
\def\IE{\relax\hbox{$\inbar\kern-.3em{\rm E}$}}
\def\IF{\relax\hbox{$\inbar\kern-.3em{\rm F}$}}
\def\IG{\relax\hbox{$\inbar\kern-.3em{\rm G}$}}
\def\IGa{\relax\hbox{${\rm I}\kern-.18em\mathcal{T}$}}
\def\IH{\relax{\rm I\kern-.18em H}}
\def\IK{\relax{\rm I\kern-.18em K}}
\def\IL{\relax{\rm I\kern-.18em L}}
\def\IP{\relax{\rm I\kern-.18em P}}
\def\IR{\relax{\rm I\kern-.18em R}}
\def\IZ{\relax{\rm Z\kern-.5em Z}}

\begin{document}
	\tikzset{
		position label/.style={
			below = 3pt,
			text height = 1.5ex,
			text depth = 1ex
		},
		brace/.style={
			decoration={brace, mirror},
			decorate
		}
	}

	\thispagestyle{empty}
	\quad
	
	\vspace{2cm}
	\begin{center}
		
		\textbf{\huge Trees with exponential height dependent weight}
		
		\vspace{1.5cm}
		
		{\Large Bergfinnur Durhuus ~~~~~~~~ Meltem Ünel}
		\vspace{0.5cm}
		
		{\it Department of Mathematical Sciences, Copenhagen University\\
			Universitetsparken 5, DK-2100 Copenhagen {\O}, Denmark}

		\vspace{0.5cm} {\sf durhuus@math.ku.dk, meltem@math.ku.dk}

		\vspace{1.5cm}
		
		{\large\textbf{Abstract}}\end{center} We consider planar rooted random trees whose distribution is even for fixed height $h$ and size $N$ and whose height dependence is of exponential form $e^{-\mu h}$. Defining the total weight for such trees of fixed size to be $Z^{(\mu)}_N$, we determine its asymptotic behaviour for large $N$, for arbitrary real values of $\mu$. Based on this we evaluate the local limit of the corresponding probability measures and find a transition at $\mu=0$ from a single spine phase to a multi-spine phase. Correspondingly, there is a transition in the volume growth rate of balls around the root as a function of radius from linear growth for $\mu<0$ to the familiar quadratic growth at $\mu=0$ and to cubic growth for $\mu> 0$.

\vspace{0.5 cm}	


	\newpage

\tableofcontents	
	
	\section{Introduction}
	
	Random trees have been at the stage of research in theoretical probability for decades, their relationship with important classes of branching processes being a strong motivating factor. In particular, the class of Bienaym{\'e}-Galton-Watson (BGW) processes and associated probability measures on trees have been intensively studied \cite{athreyaney1972branching,drmota2009random,abraham2015introduction}. In recent years, important additional motivation for studying random trees, as well as more general random graphs, originates from theoretical physics, e.g. by serving as models of statistical mechanical systems in random environments and by providing a framework for investigating non-perturbative aspects of quantum gravity \cite{ambjorn1997quantum}. Of particular relevance for two-dimensional systems are topics relating to random surfaces and random planar maps, on which significant progress has been obtained, both on the combinatorial side \cite{schaeffer1998conj,bouttier2003geodesic}, sometimes involving nontrivial correspondences between the maps in question and classes of labelled planar trees, and on the analytic side yielding constructions of interesting \emph{local limits} \cite{angel2003uniform,chassaing2006local,krikun2005local,curien2012view,menard2010two} and \emph{scaling limits}  \cite{chassaing2004random,marckert2006limit,le2007topological,le2012scaling,le2013uniqueness,miermont2013brownian,addario2017scaling}.
	
	A particularly simple and concrete correspondence between planar rooted trees and planar maps is provided by the so-called \emph{causal triangulations} of the disc \cite{malyshev2001two,durhuus2010spectral}.  The infinite volume limit of this ensemble can via the mentioned correspondence  be identified with the local limit as $N\to\infty$ of the uniform distribution of rooted planar trees of size $N$ \cite{durhuus2010spectral}, which is called the Uniform Infinite Planar Tree (UIPT) and is a special case of a more general construction of local limits of BGW measures conditioned on size  \cite{kennedy1975galton, aldous1998tree,janson2012simply}. From a physical point of view, it is natural in this context to consider the weight of a given tree $T$ to be given by 
	\bb\label{weight1} 
	w(T) = e^{-\Lambda |T|}\,,
	\ee
	where $\Lambda$ is a real constant and $|T|$ denotes the size of $T$, which also equals half the area (number of triangles) of the corresponding triangulation. We shall use the notation $g=e^{-\Lambda}$ in the main text below and will see that the total volume $X(g)$ of the measure on the space of all finite planar rooted trees defined by \eqref{weight1} is finite if and only if $g\leq \frac 14$. Written as a power series in $g$, it equals the generating function for the number $A_N$ of all rooted planar trees of size $N$, while for the particular "critical" value $g=\frac 14$ the measure equals a critical BGW measure up to a factor $2$, denoted by $\rho$ in the following (see equation \eqref{def:rho} below). When restricting the measure defined by \eqref{weight1} to the set $\cT_N$ of trees of fixed size $N$ and normalising, one evidently obtains the uniform distribution $\nu_N$ on $\cT_N$, independently of $g$, whose explicit form is given by \eqref{def:nuN} below. As already mentioned, the local limit of $\nu_N$ as $N\to\infty$ equals the UIPT, whose main features were first uncovered by Kesten \cite{kesten1986subdiffusive}. In particular, the fact that the local limiting measure is supported on trees with a single spine, i.e. a unique linear infinite subgraph emerging from the root, and that the branches attached to the spine are independently distributed according to $\rho$, are of importance for the subsequent discussion. 
	
	The form \eqref{weight1} lends itself to a natural generalisation,
	\bb\label{weight2}  
	w(T) = e^{-\Lambda |T| -\mu h(T)}\,,
	\ee
	where $\mu$ is a real constant and $h(T)$ denotes the height of $T$. Indeed, this form of weight turns out to be of relevance, when considering certain types of loop models on random causal triangulations as realised in \cite{durhuus2021critical}. It is worth noting that, for $\mu\neq 0$, weight functions of the form \eqref{weight2} do not define \emph{simply generated trees} \cite{meir1978altitude}, i.e.  $w(T)$  cannot be written as a product over vertices in $T$ of a local weight function depending only on the vertex degree. Hence, analytic tools depending on this feature, typically in the form of recursion relations, are not readily available in case $\mu\neq 0$. The goal of this paper is to present a generalisation of the local limit result for $\mu=0$ to arbitrary $\mu\in\mathbb R$ and to give a basic characterisation of the corresponding random trees, including a determination of volume growth exponents. In particular, we show that the single spine feature persists for $\mu<0$ but  the branches  become subcritical BGW trees, whereas for $\mu>0$ the spine seizes to be one-ended and becomes a random tree of its own with statistically dependent (infinite) branches and whose $n$'th generation size is Poisson distributed with mean $n\mu$. Moreover, the full local limit is in this case obtained by grafting independent critical BGW trees onto the random spine. These measures also occur as local limits of BGW trees conditioned on the asymptotic behavior of generation size, see \cite{abraham2020veryfat}.
	
	The paper is organized as follows. In subsection 2.1, we give a combinatorial definition of rooted planar trees, fix some notation, and introduce a convenient metric defining a topology and an associated Borel $\sigma$-algebra on the space of rooted planar trees, on which the probability measures in question will be defined, providing an appropriate setting for discussing local (or weak) limits. In subsection 2.2, the analytic structure of the familiar generating functions $X_m$ for the number of trees of height at most $m$ is determined and used to derive closed expressions for their Taylor coefficients $A_{m,N}$, that will be of importance for the analysis in section 4. Moreover, the probability measures  $\nu^{(\mu)} _N,\, N\in\mathbb N$, obtained by restricting the measure defined by \eqref{weight2} to trees of fixed size $N$ and normalising, are introduced. As in the case $\mu=0$, they are independent of $\Lambda$. Choosing $\Lambda=0$, the relevant normalisation factors are denoted by $Z_N^{(\mu)}$ and are referred to as \emph{finite size partition functions}. Their asymptotic behaviour for large $N$ is of crucial importance for the determination of the local limit. Section 3 covers the case $\mu<0$. A rather simple analysis of the analytic properties of the generating function $Z^{(\mu)}(g)$ for the partition functions $Z^{(\mu)}_N$ allows a  determination of their asymptotic behavior in subsection \ref{sec:3.1}. Subsequently, lower bounds on ball volumes are established in subsection \ref{sec:3.2}, that are strong enough to prove tightness of the measures $\nu^{(\mu)}_N,\, N\in\mathbb N$, as well as existence of the weak limit, including explicit expressions for measures of balls. The latter result allows an identification of the limiting measure with the local limit of the measures obtained by conditioning a subcritical BGW measure $\{h(T)\geq n\},~ n\in\mathbb N$, as detailed in Proposition \ref{corconditional}. In section 4, the case $\mu>0$ is considered. A saddle point approach is applied to determine the asymptotic behavior of $Z^{(\mu)} _N$ for large $N$ in subsection \ref{sec:4.1}, while in subsection \ref{sec:4.2} this approach is extended to obtain estimates on ball volumes that ultimately allow us to prove tightness and existence of the local limit $\nu^{(\mu)}$. The properties of $\nu^{(\mu)}$ are investigated in subsection 4.3 where we first establish a decomposition result in Theorem \ref{decomposition}, which also allows us to identify the measure governing the spine, including the mentioned Poisson distribution of its generation size, see Theorem \ref{backbonedist}. Finally, the statistical behaviour of the volume of the ball of radius $r$ around the root of a tree is investigated as a function of $r$. It is shown that its expectation value is quadratic in $r$ for the spine while it is cubic for the full measure $\nu^{(\mu)}$, see Corollaries \ref{expectedbackbone} and \ref{expected}.
	The corresponding almost sure statements are formulated in Corollary \ref{expectedbackbone} iv) and Theorem \ref{asvolgrowth}, while some technical arguments are deferred to an appendix.

	\section{Preliminaries}\label{sec:2}
	
	\subsection{Metric spaces of trees}\label{sec:2.1}
	
	Combinatorially, a planar tree $T$ is given by a sequence $D=(D_0, D_1, D_2,\dots)$ of (disjoint) ordered, finite sets whose elements are the vertices of $T$,
	$$
	V(T)= \bigcup_{r=0}^\infty D_r \,,
	$$
	\noindent and a sequence $\phi=(\phi_1,\phi_2, \phi_3,\dots)$ of order preserving maps $\phi_r: D_r\to D_{r-1}$, called \emph{parent maps}, such that the edges of $T$ are of the form $\{i,\phi_r(i)\}$, i.e.
	$$
	E(T)=\{\{i,\phi_r(i)\}\mid r\in\mathbb N,\, i\in D_r\}\,.
	$$ 
	Moreover, we assume $D_0=\{i_0\}$ and $D_1=\{i_1\}$ are one-point sets and call $i_0$ the \emph{root}  and $\{i_0,i_1\}$ the \emph{root edge} of $T$. Vertices in $\phi_r^{-1}(j)$ are referred to as \emph{off-spring of} $j\in D_{r-1}$ and they inherit an ordering from $D_{r}$. The notion of \emph{ancestor} and \emph{descendant} are defined in terms of the parent maps in the obvious way.
	Two trees $T=(D,\phi)$ and $T'=(D',\phi')$ are considered equal if there exist order preserving bijective maps $\psi_r: D_r\to D'_r$ such that $\psi_{r-1}\circ\phi_r=\phi'_r\circ\psi_r$, for each $r\in\mathbb N$.  It is easy to see that, by choosing a fixed orientation of $\mathbb R^2$ and a right-handed coordinate system, one can embed any such tree as a graph in $\mathbb R^2$ such that the vertices in $D_r$ are mapped into the horizontal line through $(0,r)$ and ordered according to their first coordinate where edges are represented by straight line segments, and two embedded trees are identical if and only if one can be mapped onto the other by an orientation preserving homeomorphism of $\mathbb R^2$. Hence, we shall frequently refer to the ordering of $D_r$ as being from left to right. 
	
	The neighbours of a vertex $i\in D_{r}, r\geq 1$,  consist of its offspring together with its parent. Hence, the \emph{degree} $\sigma_i$ of $i$ in $T$ is given by
	$$
	\sigma_i = |\phi_{r+1}^{-1}(i)|+1\,,
	$$
	where we write $|A|$ for the cardinality of a set $A$. Ordering the offspring $i_1, i_2,\dots, i_{\sigma_i-1}$ from left to right, the oriented edges $(i, i_n),\,n=1,\dots,\sigma_i-1$ will be called \emph{outgoing from} $i$. Together with $(i,\phi_r(i))$, these edges divide a small disc around $i$ into $\sigma_i$ angular sectors that will be denoted $S_{i,1},\dots,S_{i,\sigma_i}$, ordered consistently with the ordering of the offspring, and we will call them the \emph{sectors around the vertex} $i$. Combinatorially, the sector $S_{i,n}$ may be identified with the pair of oriented edges $((i,i_{n-1}),(i,i_{n}))$, for $n=1,\dots,\sigma_i$, with the convention $i_0 = i_{\sigma_i} =\phi_r(i)$.
	
	By ${\mathcal T}_N$ we shall denote the set of planar trees $T$ of size $N$, i.e. $|T|:=|E(T)|=N$ . Here $N\in{\mathbb N}$ or $N=\infty$. Thus, the set of finite trees is  
	$$
	{\mathcal T}_{\rm fin} := \bigcup\limits_{N=1}^\infty {\mathcal T}_N\,,
	$$ 
	and we define 
	$$
	\cT = \cT_{\rm fin}\cup \cT_{\infty}\,.
	$$
	
	Given $T$ as above and letting $d_T$ denote the standard graph distance on $T$,  it is clear that $D_r$ is the set of vertices at distance $r$ from the root, which we shall also call the \emph{height} of those vertices in $T$. The height $h(T)$ of a tree $T\in\cT_{\rm fin}$ is defined as  the maximal height of any vertex in $T$, or
	$$
	h(T) = {\rm max}\{r\mid D_r(T)\neq \emptyset\}\,.
	$$
	If $T\in \cT_\infty$ we set $h(T)=\infty$.
	
	With the aim of discussing local limits of sequences of measures on $\cT$ we shall equip it with a natural metric as follows.   
	Given $r\in\mathbb N$, let $B_r(T)$ denote the ball of radius $r$ around the root $i_0$ defined as the subtree of $T$ spanned by the vertices at distance at most $r$ from $i_0$, i.e.
	$$
	V(B_r(T)) = \bigcup_{s=0}^r D_s \,.
	$$
	For $T,T'\in\cT$ we then set 
	$$
	\text{dist}(T,T') = \inf\{\frac{1}{r} \mid r\in\mathbb N,\, B_r(T)=B_r(T')\} \,. 
	$$
	It is then easy to see that dist is a metric on $\cT$, in fact an ultrametric. We shall denote by $\mathcal{B}_{a}(T_0)$  the ball of radius $a>0$ around $T_0\in\cT$, 
	$$
	\cB_a(T_0) := \{ T \in \cT \mid \text{dist}(T,T_0) \leq a \} \, .
	$$
	
	The measures on $\cT$ discussed in the following are all Borel measures, i.e. they are defined on the Borel $\sigma$-algebra $\cF$ generated by the open sets. By definition, a sequence $\nu_N\,, \; N\in\mathbb N$, of probability measures converges weakly to a probability measure $\nu$ on $\cT$, if 
	$$
	\int_\cT F\,d\nu_n\;\to\; \int_\cT F\,d\nu \quad \mbox{as $N\to\infty$}
	$$
	for all real valued bounded continuous functions $F$ on $\cT$. This requirement is equivalent to the statement that
	\bb\label{convcond}
	\nu_N(B)\;\to\; \nu(B)\quad\mbox{as $N\to\infty$}
	\ee
	for any ball $B$ in $\cT$ as further detailed in the following remark, listing some basic properties of the metric space $\cT$. They are easily verifiable and will be used repeatedly in the subsequent discussion, see e.g.  \cite{chassaing2006local,durhuus2003probabilistic} for more details. 
	\begin{remark}\label{remcT}
		
		i)\;  If $T\neq T'$, then $\text{dist}(T,T') = \frac {1}{r}$, where $r\geq 1$ is the radius of the largest ball around their roots shared by $T$ and $T'$, and in this case we have 
		$$
		\cB_{\frac 1r}(T) = \cB_{\frac 1r}(T') = \cB_{\frac 1r}(T_0)\,,\quad\mbox{where $T_0=B_r(T)$ and hence $h(T_0)=r$}.
		$$ 
		Noting that $\cB_{\frac 1s}(T_0) =\{ T_0\}$ for $s>r$, it follows that any finite tree is an isolated element of $\cT$. 
		
		If integers $s, r$ fulfill $1\leq s < r$ and $T_1$ is a finite tree of height $s$, there is a unique decomposition of $\cB_{\frac 1s}(T_1) $ into disjoint balls of radius $\frac 1r$,
		$$
		\cB_{\frac 1s}(T_1) = \bigcup_{\substack{T_0\in\cT_{\rm fin},\, h(T_0)\leq r \\ B_s(T_0) = T_1}} \cB_{\frac 1r}(T_0)\,.
		$$
		
		ii)\; Any ball in $\cT$ is both open and closed and any two balls are either disjoint or one is contained in the other. Since $\cT_{\rm fin}$ is a countable dense subset of $\cT$, it follows immediately by use of Theorem 2.2 in \cite{billingsley2013convergence} that a sequence of probability measures $\nu_N,\, N\in\mathbb N$, converges weakly to a probability measure $\nu$ on $\cT$, if \eqref{convcond} holds for all balls $B$. 
		
		iii)\; For later reference we define
		$\cF_r$ to be the collection of subsets of $\cT$ that can be written as a countable union of balls of radius $\frac 1r$. It follows from i) above that $\cF_r,\, r\in\mathbb N$, is a filtration of $\sigma$-algebras consisting of sets that are both open and closed. Moreover, since $\cT$ is separable, the set algebra 
		$\cF_\infty := \cup_{r\in\mathbb N} \cF_r$  generates the Borel $\sigma$-algebra $\cF$ of $\cT$. 
		
		iv)\; $\cT$ is a complete metric space. It is not compact, but subsets the form 
		\bb \nonumber  \label{Compact}
		C = \bigcap _{r=1} ^\infty \{ T \in \mathcal T \big| ~ |D_r (T)|  \leq K_r \}\,,
		\ee
		where $(K_r)_{r\in\mathbb N}$ is any sequence of positive numbers, are compact. 
	\end{remark}
	Given $T\in\cT_\infty$, we say that a vertex $i$ of $T$ is of \emph{infinite type} if it has infinitely many descendants, and otherwise it is of \emph{finite type}. Clearly, if $i$ is of infinite type then so are all its ancestors and, since $T$ is locally finite, $i$ has at least one off-spring of infinite type. It follows that the vertices of infinite type span a subtree of $T$ with the same root and root edge and with no leaves, i.e. no vertices of degree 1. We call this subtree the  \emph{spine} or \emph{backbone} of $T$ and denote it by $\chi(T)$.  The mapping $\chi:\cT_\infty\to\cT_\infty$ will be called the \emph{spine map}, and we shall need the following fact about it.
	
	\begin{lem}
		The spine map $\chi:\cT_\infty\to\cT_\infty$ is Borel measurable.
	\end{lem}
	
	\begin{proof} Let $\chi_r:\cT_\infty\to\cT_\infty$ be defined such that $\chi_r(T)$ is obtained by deleting the leaves in $D_r(T)$ from $T$ together with the edges containing them. Then $\chi_r(T)$ is a subtree of $T$ with the same spine, since only vertices of finite type are deleted and only finitely many of them. Note also that $\chi_r$ is continuous, since $\chi_r(\cB_{\frac 1s}(T))\subseteq \cB_{\frac 1s}(\chi_r(T))$ for $s > r$. It is thus sufficient to show that 
		$$
		\chi(T)=\lim_{n\to\infty}\chi_2\circ\chi_3\circ\dots\circ\chi_n(T)\,.
		$$
		For this purpose, let $i$ be any vertex of finite type in $T$ at height $r$ and consider the finite subtree $T_0$  of $T$ spanned by the descendants of $i$ together with $i$ and $j=\phi_r(i)$, and where $j$ is the root of $T_0$. Observe that, if $k$ denotes the height of $T_0$, then $\chi_2\circ\chi_3\circ\dots\circ\chi_n(T)$ does not contain $i$ when $n\geq r+k-1$. Given $s\in\mathbb N$,  there are only finitely many vertices in $T$ of height $\leq s$, so it follows from this observation that if $n$ is large enough then $\chi_2\circ\chi_3\circ\dots\circ\chi_n(T)$ does not contain any vertices of finite type at heights $\leq s$, and hence coincides with $\chi(T)$ up to height $s$. This proves the claim. 
	\end{proof}
	
	Trees belonging to  $\chi(\cT_\infty)$ will be called \emph{spine trees} in the following and we will use the notation $\cT^s$ for this set.
	
	\begin{lem}
		The subset $\cT^s$ of $\cT_{\infty}$ is closed.
	\end{lem}
	
	\begin{proof}
		Note, that a tree $T\in\cT$ belongs to the complement of $\cT^s$ if and only if it has at least one leaf, and if $T$ has a leaf at height $r$ then the trees in $\cB_{\frac 1s}(T)$ also have a leaf at height $r$ for $s>r$. Hence the complement of $\cT^s$ is open, which proves the lemma.  
	\end{proof}
	
	We shall make use of the process of \emph{grafting} a tree $T_1\in \cT$ onto another tree $T_0\in\cT$ on several occasions in the following. Although intuitively rather clear, it will be  useful to specify this notion and the associated notation in some detail. Let $i$ be some vertex in $T_0$ at height $r_1\geq 1$ and let $1\leq n\leq\sigma_i$.  The grafted tree $T := gr(T_0;(i,n);T_1)$ is then defined by setting 
	\bb\label{graft1}
	D_r(T) = \begin{cases}D_r(T_0)\quad\mbox{ if $r\leq r_1$}\\ D_r(T_0)\cup D_{r-r_1+1}(T_1)\quad\mbox{if $r>r_1$}\end{cases}
	\ee
	with ordering induced by those of $T_0$ and $T_1$ such that, if $r>r_1$ and $j\in D_r(T_0)$, then $j$ is to the left of $D_{r-r _1+1}(T_1)$ if the ancestor of $j$ in $D_{r_1}(T_0)$ is to the left of $i$ or if it equals $i$ and $j$ is among the first $n-1$ offspring of $i$ or a descendant thereof.  Otherwise, $j$ is to the right of $D_{r-r_1+1}(T_0)$. Moreover, the parent maps of $T$ are defined in the obvious way in terms of those of $T_0$ and $T_1$, with the specification that the parent of vertices in $D_{2}(T_1)\subseteq D_{r_1+1}(T)$ is defined to be $i$.  In this way, $T_1$ can be considered as a subtree of $T$ whose root edge is identified with $\{i,\phi_{r_1}(i)\}$. Pictorially, one can think of $T$ as being obtained by identifying the outgoing  root edge of $T_1$ with $(\phi_{r_1}(i),i)$ and drawing the remaining part of $T_1$ in the $n$'th sector $S_{i,n}$ of the plane around $i$, see Fig.1. Likewise, $T_0$ is a subtree of $T$ with the same root and root edge and with vertices of identical degrees in both, except for $i$ in case $T_1$ has more than one edge. We say that $T$ is obtained by \emph{grafting $T_1$ onto $T_0$ at $(i,n)$} or in sector $S_{i,n}$. In case $i$ is a leaf, there is only one sector $S_{i,1}$ and we say that $T_1$ is grafted onto $T_0$ at $i$.
	
	It is easily seen that, for fixed $T_0$ and pairs $(i_1,n_1),\dots,(i_K,n_K)$ labelling different vertex sectors,  successive grafting of trees $T_1,\dots, T_K$ at $(i_1,n_1),\dots, (i_K,n_K)$, respectively, is well defined and independent of the order of grafting. We denote the so obtained tree by $gr(T_0;(i_1,n_1),\dots,(i_K,n_K);T_1,\dots,T_K)$. 
	
	The following result on the process of grafting will be needed. 
	
	\begin{lem}\label{lem:2.3}
		For fixed $T_0\in \cT$ and different pairs $(i_1,n_1),\dots,(i_K,n_K)$ as above, the mapping
		$$
		G: (T_1,\dots,T_K)\to gr(T_0;(i_1,n_1),\dots,(i_K,n_K);T_1,\dots,T_K)
		$$
		has the following properties.
		
		i)\; $G$ maps $\cT^K$ homeomorphically onto a closed subset of $\cT$. If $T_0\in \cT_{\rm fin}$, the image is also open. 
		
		ii)\; If $T_0$ is finite and $i_1,\dots, i_K$ denote the vertices (leaves) at maximal height $r:=h(T_0)$, then $G$ maps $\cT^K$ homeomorphically onto $\cB_{\frac 1r}(T_0)$.
	\end{lem}
	
	\begin{proof} It is clear from the definition of the grafting operation that 
		$$
		G(\cB_{\frac 1r}(T_1)\times\dots\times\cB_{\frac 1r}(T_K))\subseteq \cB_{\frac 1{r}}(G(T_1,\dots,T_K))\,.
		$$
		This shows, in particular, that $G$ is a contraction and so is continuous.  Letting $r_k$ be the height of $i_k$ in $T_0$ for $k=1,\dots,K$ and letting $\overline r$ denote the maximal of these heights, it also follows that $G(T_1,\dots,T_K)$ and $G(T'_1,\dots,T'_K)$ have a common ball around the root of radius $r\geq\overline r$ if and only if $T_k$ and $T'_k$ have a common ball of radius $r - r_k+1$, $k=1,\dots, K$. This implies that 
		$$
		\mbox{dist}(G(T_1,\dots,T_K),G(T'_1,\dots,T'_K))\geq \frac{1}{\overline r}\mbox{max}\{\mbox{dist}(T_k,T'_k)\mid k=1,\dots,K\}\,,
		$$
		from which we conclude that $G$ is injective with a continuous inverse defined on the image of $G$, which is closed by the completeness of $\cT$.
		This establishes the first part of i). 
		
		Assume now that $T_0$ is finite. The claim in ii) then follows by observing that for any tree $T\in \cB_{\frac 1r}(T_0)$ we have
		$
		T=G(T_1,\dots, T_K)\,, 
		$
		where $T_k$ is the subtree of $T$ spanned by $i_k$ and its descendants together with its root $\phi_r(i_k)$, for  $k=1,\dots, K$. 
		
		In order to verify the last statement in i), let $T=G(T_1,\dots,T_K)$. If $T$ is finite it is isolated in $\cT$ and hence an open subset. On the other hand, if  $T$ is infinite, some $T_k$ must be infinite. Considering $T_0':=B_{s}(T)$ for some integer $s>h(T_0)$ and letting $ j_1,\dots,j_M$ denote the vertices at maximal height $s$ in $T_0'$,  we have that $T$ is contained in the ball of radius $\frac{1}{s}$ around $T_0'$ and by ii) this ball is contained in the image of $G$, since any of the vertices $j_1,\dots, j_M$ is contained in some $T_k$. This completes the proof of the lemma.
	\end{proof}

	
	\begin{figure}[ht]
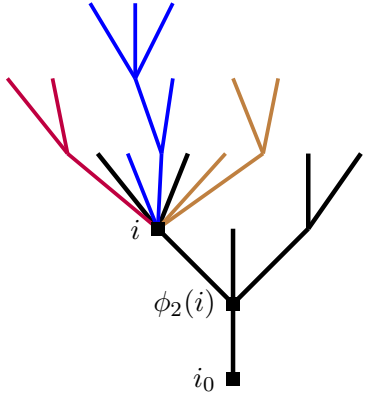

		\centering
		\begin{align*}
			\raisebox{0cm}{\tree} 
		\end{align*}
		\captionof{figure}{A finite tree $T_0$ (edges in black) with root $i_0$ and subtrees  \textcolor{purple}{$T_1$}, \textcolor{blue}{$T_2$}, \textcolor{brown}{$T_3$} grafted in the three sectors at $i$. }
		\label{figure:tree}
	\end{figure}
	
	\subsection{Generating functions}\label{sec:2.2}
	
	In this subsection we introduce the quantities needed in order to define the measures on $\cT$ whose local limits will be investigated in subsequent sections. The use of generating function techniques to deal with the combinatorial properties of those quantities is a main theme of the discussion below.
	
	The generating function for the number $A_{m,N}$ of trees in $\cT_{\rm fin}$ of height at most $m$ and size $N$ is defined by 
	\bb \nonumber \label{partfuncfinite}
	X_m(g) = \sum_{h(T)\leq m} g^{|T|} = \sum_{N=1}^\infty A_{m,N} g^N\,.
	\ee
	Similarly, we define 
	\bb\label{defX}
	X(g) =  \sum_{T\in \cT_{\rm fin}} g^{|T|} = \sum_{N=1}^\infty A_{N} g^N\,,
	\ee
	where $A_N$ is the number of trees in $\cT_{\rm fin}$ of size $N$. In order to study the convergence properties of the sums involved, it is convenient to make use of the well known recursion relation, see e.g. \cite{drmota2009random,durhuus2003probabilistic},
	\bb\label{recursion1}
	X_{m+1}(g) = \frac{g}{1-X_m(g)},\; m\geq 1, \quad X_1(g)=g\,,
	\ee
	as well as the equation satisfied by $X$,
	\bb\label{eqX}
	X(g) = \frac{g}{1-X(g)}\,.
	\ee
	The unique solution to \eqref{eqX} with $X(0) = 0$ is 
	\bb\label{solX}
	X(g) = \frac{1-\sqrt{1-4g}}{2}\,,
	\ee
	from which the values of its Taylor coefficients $A_N$ can be deduced, yielding 
	\bb\label{asympCatalan}
	A_N = C_{N-1} := \frac{(2N-2)!}{N!(N-1)!} =  \frac{1}{\sqrt\pi}N^{-\frac 32} 4^{N-1}\big(1+O(N^{-1})\big)\,,
	\ee
	where $C_N$ is known as the $N$'th Catalan number. In particular, $X$ is analytic in the disc 
	$$
	\D=\{g\in\mathbb C\mid\; |g|<\frac 14\}\,.
	$$
	Moreover, $X$  equals the limit of $X_m$ as $m\to\infty$ on $\mathbb D$, as can be seen from the explicit formula for $X_m$ stated below. This convergence will also be discussed in more detail in subsection \ref{sec:3.1}. 
	
	Equation \eqref{recursion1} can be rewritten in linear form and solved explicitly with the result  (see e.g. \cite{drmota2009random})
	\bb\label{eq:sol1}
	X_m(g) = 2g\frac{(1+\sqrt{1-4g})^m -(1-\sqrt{1-4g})^m}{(1+\sqrt{1-4g})^{m+1} - (1-\sqrt{1-4g})^{m+1}} \,,\quad m\geq 1\,.
	\ee
	Recalling that the Chebychev polynomials $U_m$ of the second kind are defined by
	\bb\label{eq:cheb}
	U_m(\cos\theta) = \frac{\sin(m+1)\theta}{\sin\theta}
	\ee
	and setting
	$$
	e^{i\theta} = \frac{1+\sqrt{1-4g}}{2\sqrt g}\,,
	$$
	one finds that
	\bb\label{eq:sol2}
	X_m(g) = \sqrt g\frac{U_{m-1}(\frac{1}{2\sqrt g})}{U_m(\frac{1}{2\sqrt g})}\,.
	\ee
	The roots $x_{m,k}$ of $U_m$ are given by 
	$$
	x_{m,k}= \cos\theta_{m,k}\,, \quad\mbox{where}\quad \theta_{m,k}=\frac{\pi k}{m+1}\,\mbox{and}\; k=1,\dots, m.
	$$ 
	Recalling that $U_m$ is of degree $m$ and has the same parity as $m$, it follows that if $m=2l$ is even we have that $g^l U_{m}(\frac{1}{2\sqrt g})$ is a polynomial in $g$ of degree $l$ with non-vanishing constant term, while if $m=2l+1$ is odd then the same holds for $\sqrt g g^l U_{m}(\frac{1}{2\sqrt g})$.  Hence, we conclude from \eqref{eq:sol1} that $X_m$ is a rational function of $g$ of the form
	\bb\label{eq:rat}
	X_m(g) = \frac{g P_m(g)}{Q_m(g)}\,,
	\ee
	where both $P_m$ and $Q_m$ are polynomials of degree $\lfloor\frac m2\rfloor$ if $m$ is odd, whereas they are of degree $\frac m2-1$ and $\frac m2$, respectively, if $m$ is even. Moreover, the roots of $Q_m$ are in both cases of the form $g_{m,k}=\frac{1}{4x_{m,k}^2}$, where $x_{m,k}$ is a nonvanishing root of $U_m$ as given above. Since $x_{m,k}=-x_{m,m+1-k}$ we obtain precisely $\lfloor\frac m2\rfloor$ different roots  $g_{m,k}$ given by 
	$$
	g_{m,k}= \frac 14\Big(1+\tan^2\frac{\pi k}{m+1}\Big)\,,\quad k=1,\dots,\big\lfloor\frac m2\big\rfloor\,.
	$$   
	These are all simple poles of $X_m$ whose corresponding residues $r_k$ we calculate below. In particular, it follows that $X_m$ is an analytic function in the disc
	$$
	\D_m := \{g\in\mathbb C\mid |g|\leq g_m\}\,,
	$$
	where 
	$$
	g_m := g_{m,1} = \frac 14\Big(1+\tan^2\frac{\pi }{m+1}\Big)
	$$
	is the radius of convergence for the power series \eqref{partfuncfinite} defining $X_m$. 
	
	Using \eqref{eq:sol2}, we obtain by differentiating the denominator on the right-hand side that 
	\bb\label{eq:res} 
	r_{m,k} = -\frac{4g_{m,k}^2U_{m-1}(\frac{1}{2\sqrt g_{m,k}})}{U'_m(\frac{1}{2\sqrt g_{m,k}})}\,.
	\ee
	From the defining relation \eqref{eq:cheb} one obtains 
	$$
	U_{m-1}(\cos\theta_k) = (-1)^{k+1}
	$$
	and
	$$
	U'_m(\cos\theta_k) = (-1)^{k+1}\frac{m+1}{\sin^2\frac{\pi k}{m+1}}\,.
	$$
	Inserting these into \eqref{eq:res} then gives 
	\bb\label{eq:res2} 
	r_{m,k} = - \frac{1}{4(m+1)}\tan^2\frac{\pi k}{m+1}\Big(1+ \tan^2\frac{\pi k}{m+1}\Big)\,,\quad k=1,\dots,\big\lfloor\frac m2\big\rfloor\,.
	\ee
	Taking into account  \eqref{eq:rat} and the relative degrees of $P_m$ and $Q_m$, we have
	\bb\label{eq:rat2}
	X_m(g) = \sum_{k=1}^{\lfloor\frac m2\rfloor} \frac{r_{m,k}}{g-g_{m,k}} + c_m +c'_mg\,,
	\ee
	where the constants $c_m$ and $c'_m$ are fixed by the requirements $X_m(0)=0$ and $X'_m(0)=1$ for $m\geq 1$. If $m$ is even, the numerator and denominator in \eqref{eq:rat} have the same degree, hence $c'_m=0$, while a simple calculation yields $c'_m = \frac{2}{m+1}$ if $m$ is odd.  
	
	By expanding the pole terms on the right-hand side of \eqref{eq:rat2} as geometric series and using \eqref{eq:res2}, we obtain the power series expansion of $X_m$, and hence its Taylor coefficients $A_{m,N}$ in the form
	\bb\label{eq:A}
	A_{m,N} = 4^N \sum_{k=1}^{\lfloor\frac m2\rfloor}\frac{1}{m+1}\tan^2\frac{\pi k}{m+1}\big(1+ \tan^2\frac{\pi k}{m+1}\big)^{-N}\,,\quad N\geq 2\,,
	\ee 
	which will constitute the basis of much of the discussion in section 4. This formula can also be found in \cite{de1972average}.
	
	\begin{remark}
		Note that, $A_{m,N}$  is constant for fixed $N$ and $m\geq N$, since the height of a tree obviously cannot exceed its size. On the other hand, the sum on the right-hand side of \eqref{eq:A} approximates the  integral
		$$
		I_N = \int_ 0^{\frac 12} \tan^2\pi x\big(1+\tan^2\pi x\big)^{-N}dx
		$$
		as $m\to\infty$. It is straight-forward to calculate $I_N$ by recursion and showing that it equals the Catalan number $C_{N-1}$ up to the factor $4^N$ in accordance with \eqref{asympCatalan}.
	\end{remark}
	
	
	The goal of the subsequent discussion is to study in some detail a one-parameter family of probability measures $\nu^{(\mu)}$ obtained as local limits of finite size measures $\nu^{(\mu)} _N$, which are uniform in size for fixed height, defined by
	\bb\label{eq:vuN}
	\nu_N^{(\mu)}(T) = \frac{e^{-\mu h(T)}}{Z^{(\mu)}_N}\,,\quad\mbox{for $T\in\cT_N$}\,,
	\ee 
	with normalisation factor (partition function) $Z^{(\mu)}_N$ given by
	\bb\label{eq:part}
	{Z^{(\mu)}_N} = \sum_{m=1}^\infty e^{-\mu m}(A_{m,N}-A_{m-1,N})\,,
	\ee 
	where $A_{0,N} =0$ by convention.
	We shall consider these as measures on the space $\cT$ and obtain the local limits as measures supported on $\cT_\infty$. For $\mu=0$, eq. \eqref{eq:part} yields the uniform probability measure on $\cT_N$ given by
	\bb\label{def:nuN}
	\nu^{(0)}_N(T)= \frac{1}{C_{N-1}}\,,\quad\mbox{for $T\in\cT_N$}\,,
	\ee
	as a consequence of \eqref{asympCatalan}. As previously mentioned, its local limit as $N\to\infty$, called the UIPT, is well studied, and a brief account of its main features is contained in the discussion in subsection \ref{sec:3.3} below, see in particular Remark \ref{muequalzero}. Our primary focus is on the cases $\mu<0$ and $\mu>0$ which, as it turns out, need to be treated by quite different techniques.
	
	\section{The case $\mu<0$}\label{sec:3}
	
	\subsection{Partition Function}\label{sec:3.1}
	
	In order to determine the asymptotic behaviour of $Z_N^{(\mu)}$ for $N$ large, as a prerequisite for establishing the existence of the local limit of $\nu_N^{(\mu)},\,N\in\mathbb N$, some more detailed information on the rate of convergence of $X_m$ towards $X$ will be needed. This is provided by the next lemma. 
	
	
	\begin{lem} \label{lemmaXm} Defining 
		\bb\label{defc}
		c(g) = \frac{g}{(1-X(g))^2}
		\ee
		we have that 
		\bb\label{cbound}
		|c(g)|\leq c(|g|) <1\quad\mbox{for}\; g\in\D\,,
		\ee
		and the following statements hold.
		\begin{itemize} 
			\item[i)]		
			\bb\label{XmXbound}
			|X_m(g)-X(g)| \leq  |g| c(|g|)^m\,,\quad m\in\mathbb N,\; g\in\D\,.
			\ee
			\item[ii)]
			The product $\prod\limits_{l=1} ^\infty \frac{1-X(g)}{1-X_l(g)} $ converges and equals an analytic function $f$ on $\D$ with no zeroes and fulfilling
			\bb\label{prodbound}
			\Big|\prod\limits_{l=1} ^m \frac{1-X(g)}{1-X_l(g)} - f(g)\Big| \leq \mbox{\rm cst}\cdot c(|g|)^m\,,\quad m\in\mathbb N,\; |g|\leq a\,,
			\ee
			for each fixed $a<\frac 14$.
		\end{itemize}
	\end{lem}
	\begin{proof} From the definition \eqref{defX} of $X$ it is clear that $|X(g)|\leq X(|g|)$, whenever the sum in \eqref{defX} is absolutely convergent, i.e. when $|g|\leq\frac 14$. Moreover, $X$ is strictly increasing on $[0,\frac 14]$ with $X(\frac 14) =\frac 12$. Hence $c$ is likewise strictly increasing with $c(0)=0$ and $c(\frac 14)=1$. From this \eqref{cbound} follows.
		
		In order to verify i), we note that 
		\bb\label{ineqX-Xm}
		|X(g)-X_m(g)| =  \big| \sum_{T\in \mathcal{T}}g^{|T|} -\sum_{\substack{T\in \mathcal{T} \\ h(T)\leq m}} g^{|T|} \big| \leq ~ X(|g|)-X_m(|g|)\,,
		\ee
		as well as the identity  
		$$
		X(g)-X_m(g) = \frac{g\big(X(g)-X_{m-1}(g)\big)}{(1-X(g))(1-X_{m-1}(g))}\,,
		$$
		which is a consequence of \eqref{recursion1} and \eqref{eqX} and by iteration gives 
		$$
		X(g)-X_m(g) = g\cdot c(g)^{m}\prod_{l=1}^{m-1} \frac{1-X(g)}{1-X_{l}(g)} \,.
		$$
		Here we note that $0\leq X_l(g) \leq X(g)<\frac 12$ for $0\leq g <\frac 14$, and hence the last product is bounded by $1$. Thus  \eqref{XmXbound} follows by use of \eqref{ineqX-Xm}. 
		
		As a consequence of \eqref{cbound} and \eqref{XmXbound} we have that $\sum_{m=1}^\infty |X_m(g)-X(g)|$ is uniformly convergent on compact subsets of $\D$. Hence, the first statements of ii) follow from Theorem 15.6 of  \cite{walter1987real}, since $X_m$ and $X$ are analytic on $\D$. Furthermore, applying standard estimates, see e.g. the proof of Theorem 15.4 of  \cite{walter1987real}, we have
		\begin{align}\nonumber
			&\Big|\prod\limits_{l=1} ^m \frac{1-X(g)}{1-X_l(g)} - f(g)\Big| = |f(g)| \Big|\prod\limits_{l=m+1} ^\infty \frac{1-X_l(g)}{1-X(g)} - 1\Big| \\ &\leq |f(g)| \Big(\prod\limits_{l=m+1} ^\infty\Big (1+\Big|\frac{X(g) -X_l(g)}{1-X(g)}\Big|\Big) - 1 \Big) \leq |f(g)| \big(e^{\sum\limits_{l=m+1}^\infty \frac{|X(g)-X_l(g)|}{|1-X(g)|}} -1\big)\,,\nonumber
		\end{align}
		from which \eqref{prodbound} follows when taking into account \eqref{XmXbound} and the continuity of $f$.
	\end{proof}

	We next consider the generating function for $Z^{(\mu)} _N$ given by
	$$
	Z^{(\mu)}(g) := \sum_{N=1} ^\infty Z^{(\mu)} _N g^N = \sum_{m=1} ^\infty e^{-\mu m} (X_m(g) -X_{m-1}(g))\,. 
	$$
	For notational convenience, we set $k:= e^{-\mu}>1$ in the remainder of this section and define 
	$$
	g_c(k)= \frac{k}{(1+k)^2}\,.
	$$  
	By direct calculation, one easily verifies that $g_c(k)\in\D$ is uniquely determined by the equation
	\bb\label{eqgc}
	k\cdot c(g_c(k)) =1\,.
	\ee
	The following theorem shows that the singularity of $Z^{(\mu)}$ closest to $0$ is shifted from $g=\frac 14$ for $Z^{(0)}=X$, to $g_c(k)<\frac 14$ for $\mu<0$, and that the singularity becomes a simple pole instead of a square root branch point.	
	
	\begin{thm}
		For fixed $k=e^{-\mu} >1$, there exists  $b >g_c(k)$ such that $Z^{(\mu)}(g)$ is analytic in
		\bb \label{puncdisk}
		\{g\in\mathbb C\mid \abs{g}< b , g \neq g_c(k) \}  \,,
		\ee
		and has a simple pole at $g_c(k)$.
		\label{thmZ}
	\end{thm}
	
	\begin{proof}
		Using \eqref{recursion1}, we have
		
		\bb
		k^m\big( X_{m+1}(g) -X_{m}(g)\big) = \frac{kg^2}{1-X(g)}  \big(kc(g)\big)^{m-1} \frac{1-X(g)}{1-X_m(g)}\prod_{l=1} ^{m-1} \Big( \frac{1-X(g)}{1-X_l(g)}\Big) ^2\,,
		\label{difference}
		\ee
		which by Lemma \ref{lemmaXm} can be written as 
		\bb\label{Zterm}
		k^m\big( X_{m+1}(g) -X_{m}(g)\big) = \frac{k(gf(g))^2}{1-X(g)}  \big(kc(g)\big)^{m-1}  + h_m(g)\,,
		\ee
		where $h_m$ is analytic in $\mathbb D$ and fulfills 
		$$
		|h_m(g)|\leq \mbox{\rm cst}\cdot \big(k c(|g|)^2\big)^m\quad\mbox{for $|g|\leq b$ and $m\in\mathbb N$}\,, 
		$$
		for any fixed  $b<\frac 14$. By \eqref{eqgc}, we have that $c(g_c(k))<1$ and hence $b>g_c(k)$ can be chosen such that $kc(|g|)^2<1$ for $|g|\leq b$.
		It then follows that  $\sum_{m=1}^\infty h_m(g)$ converges to an analytic function $h(g)$  for $|g|<b$ and hence, by summing over $m$ in \eqref{Zterm}, we conclude that
		$$
		Z^{(\mu)}(g) =  \frac{k(gf(g))^2}{\big(1-X(g)\big)\big(1-kc(g)\big)} + h(g) 
		$$
		is analytic for $|g|<b$ except at $g=g_c(k)$, which is a simple zero of the denominator $1-kc(g)$.  This completes the proof.  
	\end{proof}
	
	\begin{cor}\label{cor:Zasymp1}
		There exists $d>0$ such that
		\bb
		Z^{(\mu)}_N =  r\cdot g_c(k)^{-(N+1)} \big(1+ O(e^{-d N})\big)
		\label{Zasymp1}
		\ee
		for $N$ large, where $r$ is the residue of $-Z^{(\mu)}(g)$ at $g=g_c(k)$.
	\end{cor}
	
	\begin{proof}
		By Theorem \ref{thmZ} we may write 
		$$
		Z^{(\mu)}(g) = \frac{r}{g_c(k)-g} + \tilde h(g)\,,
		$$
		where $\tilde h$ is analytic in a disc centred at $0$ of radius $b>g_c(k)$. Expanding the pole term as a geometric series in $\frac{g}{g_c(k)}$ then yields the dominant term in \eqref{Zasymp1}, while the subdominant part arises from the Taylor coefficients of $\tilde h$.
	\end{proof}

	\subsection{Lower bounds on ball volumes and the local limit}\label{sec:3.2}

	To prove that the sequence $\nu^{(\mu)} _N, N\in\mathbb N$, given by \eqref{eq:vuN} has a weak limit on $\cT$, we proceed by first establishing lower bounds on ball volumes that will allow us to prove tightness of the measures $\nu^{(\mu)}_N$ and subsequently to identify the limit.

	\begin{lem} \label{lem:ballbound1}
		Let $k=e^{-\mu}>1$ and let $T_0\in\cT_{\rm fin}$ have height $r$ with $K$ vertices in $D_r(T_0)$. For each $M\in\mathbb N$ there exists $d>0$ such that
		\bb\label{ballbound1}
		\nu^{(\mu)} _N (\cB_{\frac 1r}(T_0)) \geq   K \cdot k^{r-1} g_c(k)^{|T_0|-K}   \Big(\sum_{S=1}^M C_{S-1}g_c(k)^{S}\Big)^{K-1} \big(1+O(e^{-d N})\big)\,.
		\ee
	\end{lem}		
	
	\begin{proof}
		Given $T_0 $ as stated, let $i_1, \dots , i_K$ denote the vertices at maximal height $r=h(T_0)$. By Lemma \ref{lem:2.3} ii), any tree $T$ in $\mathcal{B} _{\frac{1}{r}}(T_0)$ can be obtained by grafting $K$ trees $T_1, \dots, T_K \in \cT$ onto $T_0$ at $i_1,\dots, i_K$, respectively, which we shall refer to as the branches of $T$, see Fig.2. If $|T|=N$, we then have   
		$$
		N=  |T_0|+|T_1|+\dots+|T_K| -K\,,
		$$
		and hence a branch $T_j$ of maximal size among $T_1,\dots, T_K$ must fulfill
		$$
		|T_j|\geq \frac{N-|T_0|+K}{K}\,.
		$$
		Given $M\in\mathbb N$ and $j\in\{1,\dots,K\}$ and imposing the constraint $|T_i|\leq M$ for $i\neq j$, it follows that $T_j$ is the unique branch of maximal size, provided $N> K(M-1) + |T_0|$. Imposing the additional condition that $h(T_j)> M$, we obtain that $T_j$ also has maximal height among $T_1,\dots, T_K$ and that
		$$
		h(T) = h(T_j) + r-1\,.
		$$
		Hence, setting $|T_i|=N_i$, it holds that
		\bb\label{ballboundprel}
		\nu^{(\mu)} _N (\cB_{\frac 1r}(T_0)) \geq  k^{r-1} \sum_{j=1} ^K \sum_{\substack{ N_i \leq M , i\neq j\\ N_1+..+N_K= N-|T_0|+K }} (Z_N) ^{-1} Z_{M,N_{j}} \prod_{i \neq j} C_{N_i-1}  \,,
		\ee
		for $N$ large enough, where 
		$$
		Z_{M,N} := \sum\limits_{m=M+1}^\infty k^m(A_{m,N} - A_{m-1,N})\,.
		$$
		Clearly, the generating function for the $Z_{M,N}$ for fixed $M$ equals $Z(g)$ minus the function $\sum\limits_{m=1}^{M-1} k^m\big(X_m(g)-X_{m-1}(g)\big)$, which is analytic in $\D$. As a consequence, $Z_{N,M}$ has the same asymptotic form \eqref{Zasymp1} for $N$ large as $Z^{(\mu)}_N$. For any given $j$ and fixed values of $N_i,\,i\neq j$,  it follows that the corresponding term in \eqref{ballboundprel} fulfills 
		$$
		(Z_N) ^{-1} Z_{M,N_{j}} = g_c(k)^{\sum\limits_{i\neq j}N_i +|T_0|-K}\big(1+O(e^{-d N})\big)\,.
		$$
		Inserting this into \eqref{ballboundprel}, we obtain \eqref{ballbound1}.
	\end{proof}

	\begin{figure}[ht]
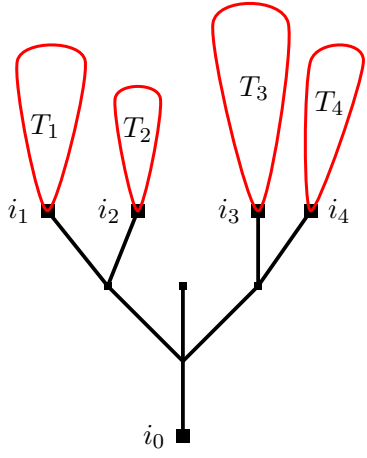

		\centering
		\begin{align*}
			\raisebox{0cm}{\treeleaf} 
		\end{align*}
		\captionof{figure}{Structure of a tree belonging to the the ball $\cB _{\frac{1}{3}}(T_0)$ around $T_0$ with root $i_0$ and leaves $(i_1,i_2,i_3,i_4)$ at height $h(T_0)=3$, obtained by grafting 4 trees $T_1,T_2,T_3,T_4 \in \cT$  onto $T_0$ at $i_1,i_2,i_3,i_4$, respectively.}
		\label{figure:treeleaf}
	\end{figure}
	
	Let us denote the the large-$N$ limit of the right-hand side of \eqref{ballbound1} by $\Lambda(T_0,M)$, i.e.
	\begin{align} \label{defLambdaM}
		\Lambda(T_0,M)  =  K \cdot k^{r-1}  g_c(k)^{|T_0| -K} \big(\sum_{S=1}^M C_{S-1} g_c(k)^S\Big) ^{K-1}\,, 
	\end{align}
	and define 
	\bb\label{defLambda}
	\Lambda(T_0) := \lim_{M\to\infty} \Lambda(T_0,M) = K \cdot k^{r-1}  g_c(k)^{|T_0| -K} X(g_c(k))^{K-1}\,. 
	\ee
	
	\begin{lem}\label{lem:addupto} For any $r\in\mathbb N$, it holds that 
		\bb\label{addupto}
		\sum_{\substack{T_0 \in \cT _{\text{\rm fin}}: h(T_0)=r}} \Lambda(T_0)  \; =\; 1\,.
		\ee
	\end{lem}
	
	\begin{proof}
		We use an inductive argument. For $r=1$ the statement trivially holds, 
		so let $r \geq 2$ be arbitrary and assume \eqref{addupto} holds for $r-1$. Recall that the factor $K$  in \eqref{defLambda} originates from summing over the positions of the long branch labelled by $j$ out of $K$ branches. This branch has a root $i$ in $D_{r-1}(T_0)$. Thus, the sum on the left-hand side of \eqref{addupto} can be rewritten as a sum over pairs $(T_0, (i,j))$ where $(i,j)$ is a marked edge at maximal height as indicated, and dropping the factor $K$. For fixed $i$, the edge $(i,j)$ divides the offspring of $i$ different from $j$ into two subsets, those to the left and those to the right of $j$, respectively, and it follows easily, for fixed $K$ and $i$, that the number of terms in the sum \eqref{addupto} equals the the number of ways of writing the number $K-1$ of edges different from $(i,j)$ as a sum of $K'+1$ non-negative integers, where $K^\prime :=|D_{r-1}(T_0)|$. Since this number equals $\binom{K+K^\prime-1}{K^\prime}$, independently of $i$, we get
		\begin{align} \nonumber
			\sum_{\substack{T_0 : B_{r-1}(T_0) = T_0 ^\prime\\ h(T_0)=r}} \Lambda(T_0) 
			&=  k^{r-1} \cdot g_c(k) ^{|T_0^\prime|}\cdot K^\prime \sum_{K \geq 1} \binom{K+K^\prime-1}{K^\prime} X(g_c(k))^{K-1} \\ \nonumber
			&= k^{r-1} \cdot K^\prime \cdot g_c(k)^{|T_0 ^\prime|} \big(1-X(g_c(k)\big) ^{-(K^\prime +1)} \\
			&= k^{r-2} \cdot K^\prime \cdot g_c(k)^{|T_0 ^\prime|-K^\prime} X(g_c(k))^{K^\prime -1}\,,  \label{identical}
		\end{align}
		where the second equality follows by using the identity
		\bb\label{combid}
		\sum_{K=R}^\infty \binom{K}{R}x^K = \frac{x^R}{(1-x)^{R+1}}\,,
		\ee
		and in the final step \eqref{eqX} and  \eqref{eqgc} were used. Since the last expression in \eqref{identical} equals $\Lambda(T_0^\prime)$, this completes the proof.
	\end{proof}	
	
	A  first consequence of Lemmas \ref{lem:ballbound1} and \ref{lem:addupto} is the tightness of the finite size measures.
	
	\begin{cor}\label{tightness0} The sequence  $(\nu_N^{(\mu)})$ of measures on $\cT$ defined by \eqref{eq:vuN} is tight for each $\mu<0$.
	\end{cor}
	
	\begin{proof}  As previously remarked, sets $C$ of the form \eqref{Compact} are compact subsets of $\cT$. Given $\epsilon>0$, it is hence sufficient to exhibit such a $C$ fulfilling $\nu^{(\mu)}_N (\cT\setminus C) <\epsilon$ for all $N$. Since 
		\bb\label{Cepsilonbound}
		\nu^{(\mu)}_N (\cT\setminus C)\;\leq\; \sum_{r=1}^\infty \nu^{(\mu)}_N (\{T\in\cT_{\rm fin} \mid |D_r(T)|> K_r\})\,,
		\ee
		if $C$ is given by \eqref{Compact}, it further suffices to show that for each $r\in\mathbb N$ and each $\delta>0$ there exists a $K>0$ such that 
		\bb\label{Ddelta}
		\nu^{(\mu)}_N (\{T\in\cT_{\rm fin} \mid |D_r(T)|> K\})\;<\; \delta\,.
		\ee
		Indeed, one can then choose $\delta$ to be $r$-dependent of the form $\delta_r=\frac{\epsilon}{2^r}$ and let  $C$ be defined in terms  of the corresponding values $K_r$ determined by \eqref{Ddelta}.  It then follows by \eqref{Cepsilonbound} that $\nu^{(\mu)}_N (\cT\setminus C) <\epsilon$.  
		
		In order to establish \eqref{Ddelta}, let $r\geq 1$ be given and choose first by Lemma \ref{lem:addupto} a finite subset $\cT _0$ of $\cT _{\text{fin}}$ consisting of trees of height $r$  such that
		\bb \label{1minusdelta}
		\sum_{T_0\in \cT_0} \Lambda(T_0) \geq 1-\delta\,,
		\ee
		and then choose by \eqref{defLambda} M large enough such that
		\bb \label{1minus2delta}
		\sum_{T_0\in \cT_0} \Lambda(T_0,M) \geq 1-2\delta\,.
		\ee
		Using \eqref{ballbound1}, we can then find $N_0$ large enough such that 
		\bb\label{1minus3delta}
		\sum_{T_0\in\cT_0}  \nu^{(\mu)} _N (\mathcal{B} _{\frac{1}{r}}(T_0) \geq 1-3\delta \quad \mbox{for $N\geq N_0$}\,.
		\ee
		Now let $K_0=\mbox{max}\{|D_r(T_0)|\mid T_0\in\cT_0\}$. Then the pairwise disjoint balls $\mathcal{B} _{\frac{1}{r}}(T_0),\, T_0\in\cT_0$, are contained in $\{T\in\cT \mid |D_r(T)|\leq  K_0\}$, and therefore  \eqref{1minus3delta} implies
		$$
		\sum_{K=1}^{K_0} \nu_N^{(\mu)}(\{T\in\cT \mid |D_r(T)| =K\}) \geq 1-3\delta \quad \mbox{for $N\geq N_0$}\,.
		$$
		Since the sets $\{T\in\cT \mid |D_r(T)| =K\},\, K\in\mathbb N$, are pairwise disjoint, it follows that
		$$
		\nu^{(\mu)}_N(\{T\in\cT \mid |D_r(T)| > K\}) \leq 3\delta\quad  \mbox{for $N\geq N_0$ and $K\geq K_0$}.
		$$
		Choosing $K> N_0$, this inequality holds for all $N$, and thus the proof is complete.
	\end{proof}
	
	\begin{thm}\label{thm:limit1}
		For each $\mu<0$ the sequence  $(\nu^{(\mu)}_N)$ is weakly convergent to a Borel probability measure $\nu^{(\mu)}$ on $\cT$, that is characterized by 
		\bb\label{limit1}
		\nu ^{(\mu)}(\cB_{\frac 1r}(T_0)) = \Lambda(T_0) = K\cdot k^{r-1}\cdot g_c(k)^{|T_0| -K} X(g_c(k))^{K-1}\,,
		\ee
		for any  tree $T_0\in\cT_{\rm fin}$ of height $r$, where $K=|D_r(T_0)|$.
	\end{thm}
	
	\begin{proof}
		Since $(\nu^{(\mu)} _N)$ is tight by the previous lemma, it has a weakly convergent subsequence $(\nu^{(\mu)} _{N_i})$ converging to a probability measure $\nu^{(\mu)} $ on $\cT$. We shall show that the limit $\nu^{(\mu)}$ is independent of the subsequence and hence that $(\nu^{(\mu)} _N)$ is convergent. Since the balls in $\cT$ have empty boundary, we  have by Theorem 2.1  in \cite{billingsley2013convergence} that $\nu^{(\mu)} _{N_i}(\cB_{\frac 1r}(T_0))$ converges to $\nu^{(\mu)}((\cB_{\frac 1r}(T_0))$ as $i\to\infty$. Using Lemma \ref{lem:ballbound1}, this implies 
		\bb \nonumber
		\nu^{(\mu)}(\mathcal{B}_{\frac{1}{r}}(T_0) )\geq \Lambda(T_0,M)
		\ee
		for any $T_0\in\cT_{\rm fin}$ of height $r \in \mathbb N$, and any $M >0$. Letting $M\to\infty$ we obtain 
		\bb \label{ineqeq}
		\nu^{(\mu)}(\mathcal{B}_{\frac{1}{r}}(T_0) )\geq \Lambda(T_0)\,.
		\ee
		Finally, using Lemma \ref{lem:addupto} and the fact that $\nu^{(\mu)}$ is a probability measure, it follows that equality holds in \eqref{ineqeq}. Since any Borel probability measure on $\cT$ is uniquely determined by its value on balls, by Theorem  2.2 in \cite{billingsley2013convergence}, this proves that the limit $\nu^{(\mu)}$ is unique. 
	\end{proof}
	
	\subsection{Properties of the local limit}\label{sec:3.3}
	In this section we discuss the measures $\nu^{(\mu)},\, \mu<0$, from the viewpoint of  branching processes, including some results on volume (or population) growth. 
	
	We first give a brief account of some aspects of Bienaymé-Galton-Watson (BGW) branching process (with one type of individual) and their local limits, while referring to, e.g., \cite{athreyaney1972branching, abraham2015introduction} for further details. Such a BGW process is defined in terms of an offspring probability distribution $p(n), n=0,1,2,\dots,$ fulfilling 
	\bb\label{offspring1}
	\sum_{n=0}^\infty p(n) =1\,.
	\ee
	The corresponding BGW tree is the probability measure $\lambda$ on $\cT$ defined by setting
	\bb\label{defGWtree1}
	\lambda(\cB_{\frac 1r}(T)) = \prod_{v\in \cup_{s=1}^{r-1}D_s(T)} p(\sigma(v)-1)\,,
	\ee
	for any $T\in\cT$. Indeed, due to \eqref{offspring1} it is straightforward to show by induction that this formula defines, for each $r\in\mathbb N$, a measure $\lambda_r$ on $\cF_r$, see Remark \ref{remcT} iii), and that they are compatible, i.e. the restriction of $\lambda_s$ to $\cF_r$ equals $\lambda_r$, if $r<s$. In particular, they define a unique finitely additive measure $\lambda_\infty$ on $\cF_\infty$. By a Kolmogoroff type of argument one can show that $\lambda_\infty$ has a unique extension $\lambda$ to a Borel measure on $\cT$.
	
	It is well known that $\lambda$ is supported on $\cT_{\rm fin}$ if and only if  $p$ is subcritical or critical, i.e. if the average offspring $m$ satisfies 
	$$
	m :=\sum_{n=0}^\infty np(n) \leq 1\,.
	$$
	Moreover, in this case it is a fundamental result of Kesten \cite{kesten1986subdiffusive} that if $\lambda_{\geq N}$ denotes $\lambda$ conditioned on $\{T\mid h(T)\geq N\}$, then the local limit of $\lambda_{\geq N}$ as $N\to\infty$
	exists and equals a BGW tree $\hat\lambda$ with two types of individuals called \emph{special} and \emph{normal}, respectively,  and whose offspring probabilities are restricted such that normal individuals can have $n=0, 1,2,\dots$ normal offspring with probability $p(n)$, but no special offspring, while special individuals can have exactly one special off-spring and in addition a number of $n\geq 0$  normal offspring with probability $p^*(n)$, given by 
	$$
	p^*(n) =  \frac{(n+1)p(n+1)}{m}\,. 
	$$
	More precisely, $\hat\lambda$ is supported on $\cT_\infty$ and its value on any ball  $\cB_{\frac 1r}(T_0)$, where $T_0\in\cT_{\rm fin}$ has height $r$. First, name each vertex of $T_0$ different from the root $i_0$ as either special or normal, such that the vertex $i_1$ next to the root in special and such that the restrictions on the offspring probabilities mentioned are respected, i.e. a special vertex has exactly one special offspring and a normal vertex has no special offspring. Clearly, the special vertices will form a path of length $r-1$ from $i_1$ to a vertex $i$ at height $r$. We call this path $\omega(i)$ since it is uniquely determined by $i$. Then associate a factor $p(\sigma(v)-1)$ with any normal vertex $v$ at height $\leq r-1$ and a  factor  $p^s(\sigma(v)-2)$ if $v$ is special, where  
	\bb\label{defps}
	p^s(n) := \frac{p^*(n)}{n+1} = \frac{p(n+1)}{m}\,,\quad n=0,1,2,\dots,
	\ee
	can be interpreted as the probability for a special individual to have $n+1$ offspring among which the special one is the $i$'th from the left, independent of $i=1,\dots, n+1$. Finally, we sum the resulting product over paths $\omega$, thus defining  
	\bb\label{defGWtree2}
	\hat\lambda(\cB_{\frac 1r}(T_0)) = \sum_{i\in D_r(T_0)}\prod_{v\in\cup_{s=1}^{r-1} D_s(T_0)\setminus V(\omega(i))} p(\sigma(v)-1)\prod_{v\in V(\omega(i))\setminus\{i\}} p^s(\sigma(v)-2)\,,
	\ee
	which characterises the measure $\hat\lambda$ on $\cT_\infty$ uniquely. It is known \cite{kesten1986subdiffusive} that $\hat\lambda$ is supported on the subset of $\cT_\infty$ consisting of trees with a single spine, i.e. the spine map $\chi$ is in this case a.s. constant and equals the infinite linear path emerging from the root (spanned by the special vertices). Due to the multiplicative structure of \eqref{defGWtree2}, the branches grafted onto the vertices next to the spine, to the left and right, are independent with identical distribution equal to the subcritical BGW tree with offspring probability $p$, while the degrees $\sigma(v)$ of the spine vertices $v$ are likewise independently and identically distributed according to $p^*(\sigma(v)-2), \sigma(v)\geq 2$. In the following we shall denote single spine trees distributed according to $\hat\lambda$ by $\hat T$.  
	
	The interpretation of $\nu^{(\mu)}$ as a BGW tree is given by the following proposition.
	
	\begin{prop} \label{corconditional}
		The measure $\nu^{(\mu)}$, for $k=e^{-\mu}>1$, equals the BGW measure $\hat\lambda$ defined above corresponding to the subcritical BGW tree with offspring probabilities given by
		\bb \label{branchingprob}
		p(n) = X(g_c(k))^n \big(1-X(g_c(k))\big) ~~,~~ n=0,1,2,\dots .
		\ee
	\end{prop}
	\begin{proof}
		Since $X(g_c(k))<\frac 12$, it is clear that  $p$ given by \eqref{branchingprob} defines a probability distribution with mean 
		\bb\label{mean}
		m=\frac{X(g_c(k))}{1-X(g_c(k))} < 1\,,
		\ee
		and hence defines a subcritical BGW process. Moreover, by \eqref{eqgc} the corresponding distribution $p^s$ has the form 
		$$
		p^s(n) = kg_c(k)\cdot X(g_c(k))^{n}\,.
		$$ 
		Using these in \eqref{defGWtree2}, it is seen for a given $T_0$ of height $r$ and with $|D_r(T_0)|=K$ that all summands have the same value equal to 
		\begin{align}\nonumber
			&~~~\big(kg_c(k)\big)^{r-1}X(g_c(k))^{|E(T_0)|-r}\big(1-X(g_c(k))\big)^{|E(T_0)|-K-r+1}\\\nonumber& = \big(kg_c(k)\big)^{r-1}g_c(k)^{|E(T_0)|-K-r+1}X(g_c(k))^{K-1}\\& = k^{r-1}\cdot g_c(k)^{|E(T_0)|-K}X(g_c(k))^{K-1}\,, \nonumber
		\end{align}
		where  \eqref{eqX} has been used in the second step. 
		Upon multiplication by $K$, the last expression is seen to coincide with \eqref{limit1}, which completes the proof of the proposition.	\end{proof}
	
	Letting $\mathbb E_\mu$ denote the expectation w.r.t. $\nu ^{(\mu)} $, the following result on the average growth of balls around the root of $\hat T$ as a function of radius is an easy consequence of familiar results about BGW processes.
	
	\begin{cor}\label{growth1} For $k=e^{-\mu}>1$, the following asymptotic relations hold, 
		\begin{eqnarray}
			\mathbb E_{\mu}(|D_r|) &=&  \frac{1+m}{1-m} + O(m^{r-1})\,,\label{growthD1av}\\
			\mathbb E_{\mu}(|B_r|) &=&  \frac{1+m}{1-m}\cdot r+ O(1)\label{growthB1av}\,,
		\end{eqnarray}
		where $m$ is given by \eqref{mean}.
	\end{cor}
	\begin{proof}  Let $T_s^L$ and $T_s^R$ be the branches of $\hat T$ grafted to the left and right, respectively, at the spine vertex at distance $s$ from the root. It follows from Proposition \ref{corconditional} and the preceding remarks that these are i.i.d. according to the subcritical BGW tree with offspring distribution given by \eqref{branchingprob}. In particular, we have (see e.g. Ch. 1 of \cite{athreyaney1972branching}) 
		$$
		\mathbb E_\mu(|D_r(T_s^L)|) = \mathbb E_\mu(|D_r(T_s^R)|) = m^{r-1}\,,\quad r\geq 1\,.
		$$
		Since 
		\bb\label{decompDr1}
		|D_r(\hat T)| = 1 + \sum_{s=1}^{r-1}\Big( |D_{r-s+1}(T_s^L)| + |D_{r-s+1}(T_s^R)|\Big)\,,
		\ee
		it follows  that 
		$$
		\mathbb E_\mu(|D_r(\hat T)|) = 1+ 2m\frac{1-m^{r-1}}{1-m}\,,
		$$
		from which the first relation follows. Using 
		\bb\label{DB}
		|B_r(\hat T)| = \sum_{s=1}^r |D_s(\hat T)|\,,
		\ee
		eq. \eqref{growthB1av} follows from \eqref{growthD1av}. 	\end{proof}
	
	A result similar to the following on a.s. asymptotic growth of individual trees should be available in the literature. For the sake of completeness, a proof is included in the appendix. 
	
	\begin{prop}\label{growth2}
		Let $k=e^{-\mu}>1$. There exist  constants $C_1, C_2>0$ and for $\nu^{(\mu)}$-a.e.  $\hat T$ a number $r_0(\hat T) \in\mathbb N$, such that 
		\begin{align}\label{eq:growth2}
			1&\leq |D_r(\hat T)| \leq C_1 \cdot \ln r\quad \mbox{and}\quad 
			r\leq |B_r(\hat T)| \leq C_2 \cdot r\ln r
		\end{align}
		for all $r\geq r_0(\hat T)$.
	\end{prop}

	\begin{remark}
		The volume growth exponent $d_h$ of a tree $T\in\cT_\infty$, defined by
		\bb\label{defHausdorff}
		d_h := \lim_{r\to\infty} \frac{\ln |B_r(T)|}{\ln r}\,,
		\ee
		is commonly referred to as the Hausdorff dimension of $T$, provided the limit exists. Proposition \ref{growth2} shows that $d_h=1$ $\nu^{(\mu)}$-a.s. for $\mu<0$.
	\end{remark}
	
	\begin{remark}\label{muequalzero}
		The local limit result of Kesten described above still holds if $\lambda_N$ is replaced by $\lambda$ conditioned on $\{h(T)=N\}$, and in case $\lambda$ is critical, i.e. if $m=1$,  the same conclusion holds if conditioning on $\{|T|=N\}$ is used instead, see \cite{abraham2015introduction} for a more detailed account. For the critical BGW process with offspring distribution
		\bb\label{offspring2}
		p_n = 2^{-(n+1)}, ~n=0,1,2,\dots, 
		\ee
		we shall denote by  $\rho$ the corresponding BGW measure on $\cT$ as defined by \eqref{defGWtree1}, which in this case reads 
		\bb\label{rhoball}
		\rho(\cB_{\frac 1r}(T_0)) = 4^{-|T_0|}\,2^{K+1}\,,
		\ee
		for any given tree $T_0\in\cT_{\rm fin}$ of height $r$ with $K$ vertices (leaves) at height $r$. As mentioned earlier, $\rho$ is supported on $\cT_{\rm fin}$ and is alternatively given by
		\bb \label{def:rho}
		\rho(T)=2\cdot 4^{-|T|} ~, ~T \in \cT_{\text{fin}} \,. 
		\ee
		By conditioning on size, it follows that 
		$$
		\rho(\cdot\mid |T|=N) = \nu^{(0)}_N\,,
		$$
		and hence the local limit $\nu^{(0)}$ in this case is supported on single spine trees $\hat T$ as before, but with branches grafted onto the spine vertices that are i.i.d. according to $\rho$ \cite{durhuus2003probabilistic}. Since $\rho$ is associated with a  critical BGW process, the volume growth exponent in this case turns out to be $d_h=2$. More precisely, it  is proven in \cite{durhuus2010spectral}, see also \cite{barlow2006random} for related results, that there exist constants $C'_1, C'_2>0$ and $r_0(\hat T)\in\mathbb N$  for $\nu^{(0)}$-a.e. $\hat T$ , such that 
		$$
		C'_1\cdot  r^2(\ln r)^{-2} \leq |B_r(\hat T)| \leq C'_2\cdot r^2\ln r\,,\quad \mbox{for $r\geq r_0(\hat T)$}.
		$$
	\end{remark}
	
	\section{The case $\mu>0$}\label{sec:4}
	
	\subsection{The partition function}\label{sec:4.1}
	
	In this subsection we determine the asymptotic behaviour of the partition functions $Z_N^{(\mu)}$ for large $N$, in case $\mu>0$, which has previously been obtained in \cite{guttmann2015analysis} (see section 7) in a slightly different context. For later purposes we provide a proof of this result. Note, however, that the order of our remainder term deviates from the one stated in \cite{guttmann2015analysis}.
	
	\begin{prop}\label{thm:4.1} For each $\mu>0$ it holds for any $\delta\in ]0,\frac 16[$ that
		\bb\label{eq:asympZ}
		Z^{(\mu)}_N = (e^\mu -1) \sqrt{\frac{\pi}{B}}\frac{\mu}{2} e^{-A N^{\frac 13}}N^{-\frac 56}4^N\Big(1+O(N^{-\delta})\Big) \,,
		\ee
		for $N$ large, where
		$$
		A =  3\Big(\frac{\pi\mu}{2}\Big)^{\frac 23} \quad\mbox{and}\quad B = 3\Big(\frac{\mu^2}{4\pi}\Big)^{\frac 23} \,.
		$$
	\end{prop}
	
	\begin{proof} Rewriting 
		\bb\label{eq:ZW}
		Z^{(\mu)}_N = (1-e^{-\mu})W_N + C_{N-1}e^{-\mu(N+1)}\,,
		\ee
		where 
		\bb\label{defWN}
		W_N = \sum_{m=1}^N e^{-\mu m} A_{m,N} \,,
		\ee
		and noting that the last term in \eqref{eq:ZW} is exponentially suppressed compared to the right-hand side of \eqref{eq:asympZ}, as a consequence of \eqref{asympCatalan}, we need to show that the asymptotic behaviour of $W_N$ is given by \eqref{eq:asympZ} with the factor $(1-e^\mu)$ replaced by $e^\mu$.
		In order to achieve this is we use a saddlepoint argument the details of which are as follows. 
		
		We first consider the contribution $\tilde W_N$ to $W_N$ obtained by  retaining only  the first term in the sum \eqref{eq:A} defining  $A_{m,N}$, and write it as 
		\bb\label{deftildeWN}
		\tilde W_N := 4^N\sum_{m=2}^N \frac{e^\mu}{m+1}\tan^2\frac{\pi}{m+1}\; e^{-f_N(m+1)}\,,
		\ee
		where 
		\bb\label{eq:deffN}  
		f_N(t) = \mu t + N\ln\Big( 1+\tan^2\frac{\pi}{t}\Big)\,,\quad t>2\,.
		\ee
		Since 
		$$
		f_N'(t) = \mu - 2\pi N\frac{1}{t^2}\tan\frac{\pi}{t}
		$$
		is a strictly increasing function of $t$ that maps  $]-\infty,2[$ onto $]-\infty,\mu[$ it follows that $f_N$ has a unique minimum $t_0$ determined as a smooth function of $\frac{\mu}{N}$ by 
		$$
		\frac{\mu}{N} = 2\pi \frac{1}{t_0^2}\tan\frac{\pi}{t_0}\,.
		$$
		Clearly, $t_0\to\infty$ as $N\to\infty$ for fixed $\mu$ and by Taylor expanding $\tan\frac{\pi}{t_0}$ on the right-hand side, we get
		$$
		\frac{\mu}{2\pi^2 N} = \frac{1}{t_0^3}\Big(1+\frac{\pi^2}{3 t_0^2} + O\big(\frac{1}{t_0^4}\big)\Big)
		$$
		which gives 
		\bb\label{eq:minasymp}
		t_0 = \Big(\frac{2\pi^2 N}{\mu}\Big)^{\frac 13} + O(N^{-\frac 13})\,.
		\ee
		Here, the leading term on the right-hand side is the unique minimum of the function
		$$
		g_N(t) = \mu t + N\frac{\pi^2}{t^2}\,,\quad t>0\,,
		$$ 
		which is obtained by Taylor expanding the $\ln$-term  in \eqref{eq:deffN} to first order in $\frac 1t$, and we thus have 
		\bb\label{eq:fg}
		f_N(t) = g_N(t) + NO\big(\frac{1}{t^4}\big)
		\ee
		for $t$ large. Let  $x_0$ denote the minimum of $g_1$, i.e.
		\bb\label{eq:x0}
		x_0=\Big(\frac{2\pi^2 }{\mu}\Big)^{\frac 13}\,,
		\ee
		and observe that 
		\bb\label{eq:gscale}
		g_N(t) = N^{\frac 13} g_1(N^{-\frac 13}t)\,.
		\ee 
		Given $\delta\in]0,\frac 16[$, we now first estimate the contribution to the sum in \eqref{deftildeWN} from terms with $m+1\geq t_0+N^{\frac 13-\delta}$. For such $m$ we have
		$$
		f_N(m+1)\geq f_N(t_0+N^{\frac 13-\delta}) = N^{\frac 13}g_1\big(x_0 + N^{-\delta} + O(N^{-\frac 23}) \big)+ NO(N^{-\frac 43})\,,
		$$
		by use of \eqref{eq:minasymp}, \eqref{eq:fg} and \eqref{eq:gscale}. Inserting the Taylor expansion
		\bb\label{taylorg1}
		g_1(x) = A + B(x-x_0)^2 + O(|x-x_0|^3)\,,
		\ee
		where
		$$
		A = g_1(x_0) = 3\Big(\frac{\pi\mu}{2}\Big)^{\frac 23}\quad\mbox{and}\quad B = \frac 12 g_1''(x_0) = 3\frac{\pi^2}{x_0^4} = 3\Big(\frac{\mu^2}{4\pi}\Big)^{\frac 23}\,,
		$$
		we thus obtain 
		$$
		f_{N+1}(m+1) \geq AN^{\frac 13}  + B N^{\frac 13-2\delta} + O(N^{-\frac 13}) + O(N^{\frac{1}{3}-3\delta})\,.
		$$
		Assuming $\delta > \frac{1}{9}$, it follows that
		\bb\label{bound1}
		\sum_{t_0 + N^{\frac 13 -\delta}\leq m+1\leq N+1}\frac{1}{m+1}\tan^2\frac{\pi}{m+1}\,e^{-f_N(m+1)} \leq {\rm cst}\cdot Ne^{-AN^{\frac 13}  - B N^{\frac 13-2\delta}}\,.
		\ee
		By the same arguments, this bound holds for the sum over $m+1\leq t_0-N^{\frac 13-\delta}$ as well.
		
		In the remaining range $t_0-N^{\frac 13-\delta}<m+1< t_0+N^{\frac 13 -\delta}$ we set 
		\bb\label{eq:xm}
		x_m = N^{-\frac 13}(m+1)\,, 
		\ee
		such that $x_m = x_0 + O(N^{-\delta})$ by \eqref{eq:minasymp}. Using \eqref{eq:fg} and \eqref{taylorg1}, we then get 
		$$
		f_N(m+1) = N^{\frac 13}g_1(x_m) + NO(N^{-\frac 43}) = N^{\frac 13}\big(A+B(x_m-x_0)^2 + O(N^{-3\delta})\big)+ O(N^{-\frac 13})\,. 
		$$
		Since $\delta>\frac 19$ this implies 
		$$
		e^{-f_N(m+1)}  = e^{-AN^{\frac 13}} e^{-B(N^{\frac 16}(x_m-x_0))^2}\big(1+O(N^{\frac 13-3\delta})\big)\,.
		$$
		Now, note that for $m$ in the range under consideration it holds that
		\bb\label{tanmscale}
		\frac{1}{m+1}\tan^2\frac{\pi}{m+1} = \frac 1N \frac{\pi^2}{x_0^3}\big(1+O(N^{-\delta})\big)\,,
		\ee
		such that we obtain
		\begin{align}\label{bound2}
			& \sum_{|m+1 - t_0|< N^{\frac 13 -\delta}}\frac{1}{m+1}\tan^2\frac{\pi}{m+1}\,e^{-f_N(m+1)}\nonumber\\
			& =   \frac{\pi^2}{x_0^3} \frac 1N  e^{-AN^{\frac 13}} \sum_{| m+1 - t_0 |<  N^{\frac 13 -\delta}}e^{ - B (N^{\frac 16} (x_m-x_0))^2}(1+O(N^{\frac13-3\delta}))\,,
		\end{align}
		since  $\delta>3\delta -\frac 13$.
		Recalling \eqref{eq:minasymp}, we see that the numbers $N^{\frac 16}(x_m-x_0)$ form a division af an interval $[a_N,b_N]$ into subintervals of length $N^{-\frac 16}$, where $a_N = -N^{\frac 16-\delta} + O(N^{-\half})$ and  $b_N = N^{\frac 16-\delta} + O(N^{-\half})$.  Therefore, disregarding the $O(N^{\frac13-3\delta})$-term, the last sum multiplied by $N^{-\frac 16}$ approximates the integral of $e^{-Bx^2}$ over $[a_N,b_N]$,
		whose value approaches $\sqrt{\frac{\pi}{B}}$ for $N$ large. Indeed, by a standard argument, we have
		$$
		\Big|\, N^{-\frac 16}\sum_{| m+1 - t_0 |<  N^{\frac 13 -\delta}}e^{ - B (N^{\frac 16} (x_m-x_0))^2} - \sqrt{\frac{\pi}{B}}\,\Big|\, \leq \, O(N^{-\delta})\,,
		$$
		and so  \eqref{bound2} implies that 
		\bb\label{bound3}
		\sum_{|m+1 - t_0|< N^{\frac 13 -\delta}}\frac{1}{m+1}\tan^2\frac{\pi}{m+1}\,e^{-f_N(m+1)}
		= \sqrt{\frac{\pi}{B}} \frac{\pi^2}{x_0^3} e^{-AN^{\frac 13}} N^{-\frac 56} \big(1+O(N^{\frac13-3\delta})\big)\,.
		\ee 
		Taking into account \eqref{bound1} and \eqref{bound3} we have shown that 
		\bb\label{tildeWasymp}
		\tilde W_N = e^\mu  \sqrt{\frac{\pi}{B}} \frac{\pi^2}{x_0^3} e^{-AN^{\frac 13}} N^{-\frac 56} 4^N \big(1+O(N^{\frac13-3\delta})\big)\,.
		\ee
		Noting that 
		\bb\label{x0mu}
		\frac{\pi^2}{x_0^3} =\frac {\mu}{2}\,,
		\ee
		and that the function $\delta\to3\delta-\frac 13$ maps $]\frac 19,\frac 16[$ onto $]0,\frac 16[$, this  is precisely of the form claimed for $W_N$. 
		
		It remains to show that  $W_N-\tilde W_N = o(\tilde W_N)$ as $N\to\infty$. For this purpose, note  that 
		$$
		\sum_{2\leq k\leq \lfloor\frac m2\rfloor}\frac{1}{m+1}\tan^2\frac{\pi k}{m+1}\Big(1+\tan^2\frac{\pi k}{m+1}\Big)^{-N} \leq \Big(1+\tan^2\frac{2\pi }{m+1}\Big)^{-(N-1)}\,,
		$$ 
		such that, by setting 
		$$
		U_N = 4^N \sum_{m=4}^N \Big(1+\tan^2\frac{2\pi }{m+1}\Big)^{-(N-1)}\,,
		$$
		we have 
		\bb\label{ineqWU}
		\tilde W_N \leq W_N \leq \tilde W_N + U_N\,.
		\ee
		Evidently, we can apply the same procedure as above to estimate $U_N$. By inspection, it is easily verified that this yields the bound
		\bb\label{bound4}
		U_N \leq \mbox{cst}\cdot e^{-h(y_0)N^{\frac 13}}N^{\frac 16} 4^N\,,
		\ee
		where 
		$$
		h(y) = \mu y + \frac{(2\pi)^2}{y^2}\,,\quad y>0\,,
		$$
		and $y_0$ is the unique minimum of $h$ determined by
		$$
		\mu = 2\frac{(2\pi)^2}{y_0^3}\,.
		$$ 
		One finds that
		$$
		h(y_0) = 3(\pi\mu)^{\frac 23} > A\,,
		$$
		and hence it follows from \eqref{tildeWasymp} and \eqref{bound4} that 
		$$
		U_N \leq  \tilde W_N\cdot e^{-cN^{\frac 13}}
		$$
		for $N$ large, where $c$ is a positive constant. Using \eqref{ineqWU}, we conclude that $\tilde W_N$ and $W_N$ have the same asymptotic behaviour given by \eqref{tildeWasymp}, and this completes the proof of the theorem.
	\end{proof}
	
	\subsection{Lower bounds on ball volumes and the local limit}\label{sec:4.2}
	We next establish lower bounds on the $\nu^{(\mu)}_N$-volume of balls that will allow us to prove tightness of the sequence $(\nu^{(\mu)}_N)$, and also to show weak convergence.
	
	\begin{lem}\label{lem:4.2} Assume $\mu>0$ and that  $T_0\in \cT_{\rm fin}$ has height $r$, and set $K = |D_r(T_0)|$.  
		Given  $0<\epsilon<\frac 1K$ and $M\in\mathbb N$, it holds for any $\delta\in ]0,\frac 16[$ that 
		\begin{align}\nonumber
			\nu^{(\mu)}_N(\cB_{\frac 1r}(T_0)) ~\geq~ &\frac{e^{-\mu(r-1)}}{4^{|T_0|-K}} \sum_{R=1}^{K} \binom{K}{R}\frac{1}{(R-1)!}\Big(\frac{\mu}{2}\Big)^{R-1}\\ & \Big(\sum_{S=1}^M C_{S-1}4^{-S}\Big)^{K-R}(1-\epsilon K)^K(1+O(N^{-\delta})) \label{est1}
		\end{align}
		for $N$ large.  
	\end{lem}
	\begin{proof} 
		Below we use $O$ to indicate a generic $O$-function satisfying $|O(x)|\leq {\rm cst}\cdot |x|$ for $x$ small enough, where the constant may depend on $K, M$ and $\epsilon$. 
		
		By Lemma \ref{lem:2.3}, the elements $T$ of $\cB_{\frac 1r}(T_0)$ are obtained by grafting $K$ trees $T_1\dots, T_K$ onto $T_0$ at the $K$ vertices of maximal height, see Fig. 2. We call $T_1,\dots, T_K$ the branches of $T$. For $T\in\cT_N$ we then have 
		$$
		N= \sum_{i=1}^K |T_i| + |T_0| - K\,,
		$$
		and we shall call the branch $T_i$ \emph{small} if $|T_i|\leq M$, while we call it \emph{large} if $|T_i| > \epsilon N$.  Note that, if $N>\frac{M}{\epsilon}$, no $T_i$ can be both small and large, and if $N>K\cdot M +|T_0|-K$ there must be at least one large branch $T_i$. Assuming $N$ fulfills these inequalities in the following, let $\Omega_{T_0,R}\,, 1\leq R\leq K$, denote the subset of $\cB_{\frac 1r}(T_0)$ consisting of trees $T$ of size $N$ whose large branches are precisely $T_1,\dots, T_R$. Since $\nu^{(\mu)}_N$ restricted to $\cB_{\frac 1r}(T_0)$ is invariant under permutation of the branches we have
		\bb\label{Rsum}
		\nu^{(\mu)}_N(\cB_{\frac 1r}(T_0) = \sum_{R=1}^K \binom{K}{R} \nu^{(\mu)}_N(\Omega_{T_0,R})\,.
		\ee
		Hence, we proceed to estimate   $\nu^{(\mu)}_N(\Omega_{T_0,R})$. 
		
		Imposing the restriction  $h(T)\geq M+r-1$ on $T\in\Omega_{T_0,R}$ ensures that the highest branch must be among the first $R$ branches. Denoting the smallest $i$ such that $T_i$ has maximal height by $j$ and $|T_i|$ by $N_i$, we have
		\begin{align}\nonumber
			&\nu^{(\mu)}_N(\Omega_{T_0,R})\,  \geq\, \big(Z^{(\mu)}_N\big)^{-1} \sum_{m=M}^\infty e^{-\mu (m+r-1)} \sum_{N_{R+1},\dots, N_K\leq M} C_{N_{R+1}-1}\dots C_{N_K-1}\\ \nonumber
			& \sum_{j=1}^R \sum_{\substack{N_1+\dots +N_R \\ = N-N_{R+1}-\dots -N_K+K-|T_0|\\ N_1,\dots, N_R > \epsilon N}} A_{m-1,N_1}\dots A_{m-1,N_{j-1}}(A_{m,N_j}-A_{m-1,N_j}) A_{m,N_{j+1}}\dots A_{m,N_R}\\
			&~~~~~~~~~~~~~~ = \frac{e^{-\mu(r-1)}}{ Z^{(\mu)}_N} \sum_{m=M}^\infty \sum_{N_{R+1},\dots, N_K\leq M} C_{N_{R+1}-1}\dots C_{N_K-1}\nonumber\\
			& ~~~~~~~~~~~~~~\sum_{\substack{N_1+\dots +N_R \\ = N-N_{R+1}-\dots -N_K+K-|T_0|\\ N_1,\dots, N_R > \epsilon N}} e^{-\mu m}(A_{m,N_1}\dots A_{m,N_{R}} -A_{m-1,N_1}\dots A_{m-1,N_R})\, . \label{bound0}
		\end{align}  
		Now, consider the last sum for fixed values of $N_{R+1},\dots ,N_K$ and set 
		$$
		N' = N_{R+1}+\dots+N_K +|T_0|-K\,.
		$$
		With notation as in subsection \ref{sec:4.1}, a lower bound on this sum
		is obtained by restricting the sum to values of $m$ fulfilling $|x_m-x_0|\leq N^{\delta}$, where $0<\delta<\frac 16$.  Hence, we proceed to further estimate this lower bound
		\begin{align}\label{defVN}
			&V_N := \sum_{|x_m-x_0|\leq N^{-\delta}} \sum_{\substack{N_1+\dots +N_R =N-N'\\ N_1,\dots, N_R > \epsilon N}} e^{-\mu m}\prod_{s=1}^R A_{m,N_s}\Big(1 - \prod_{t=1}^R\frac{A_{m-1,N_t}}{A_{m,N_t}}\Big)\,.
		\end{align} 
		Recalling the definition \eqref{eq:A} of $A_{m,N}$, we write
		\bb\label{eq:A1}
		A_{m,N_s} =  4^{N_s}\frac{1}{m+1}\tan^2\frac{\pi}{m+1}\big(1+ \tan^2\frac{\pi}{m+1}\big)^{-N_s} \Big(1+\sum_{k=2}^{\lfloor\frac m2\rfloor}\frac{\tan^2\frac{\pi k}{m+1}}{\tan^2\frac{\pi}{m+1}}\Big(\frac{1+\tan^2\frac{\pi}{m+1}}{1+\tan^2\frac{\pi k}{m+1}}\Big)^{N_s}\Big)\,.
		\ee
		Using that  
		$$
		\frac{1+\tan^2\frac{\pi}{m+1}}{1+\tan^2\frac{\pi k}{m+1}} \leq \frac{1+\tan^2\frac{\pi}{m+1}}{1+\tan^2\frac{2\pi }{m+1}} = e^{-\frac{3\pi^2}{(m+1)^2}\big(1+O(\frac{1}{m^2})\big)}
		$$
		for $2\leq k\leq \lfloor\frac m2\rfloor$, the sum in \eqref{eq:A1} is bounded from above by
		$$
		(m+1)^3 e^{-\frac{3\pi^2(N_s-1)}{(m+1)^2}\big(1+O(\frac{1}{m^2})\big)}\,. 
		$$
		Since $m = x_0 N^{\frac 13} + O(N^{\frac 13-\delta})$ for the range of $m$ considered and $N_s > \epsilon N$, we get that 
		\bb\label{eq:A2}
		A_{m,N_s} =  4^{N_s}\frac{1}{m+1}\tan^2\frac{\pi}{m+1}\big(1+ \tan^2\frac{\pi}{m+1}\big)^{-N_s}\big(1+ O(e^{-\epsilon c N^{\frac 13}})\big)\,,
		\ee
		where $c>0$ is a numerical constant. In particular, the first product in \eqref{defVN} takes the form 
		\bb\label{eq:prod1}
		\prod_{s=1}^R A_{m,N_s} =  4^{N-N'}\Big(\frac{1}{m+1}\tan^2\frac{\pi}{m+1}\Big)^R\big(1+ \tan^2\frac{\pi}{m+1}\big)^{-(N-N')}\big(1+ O(e^{-\epsilon c N^{\frac 13}})\big)^R\,.
		\ee
		To deal with the second product, we use
		\begin{align}\nonumber
			\frac{1+\tan^2\frac{\pi}{m+1}}{1+\tan^2\frac{\pi}{m}}& = e^{\frac{\pi^2}{(m+1)^2}-\frac{\pi^2}{m^2} + O(\frac{1}{m^4})}
			= e^{-\frac{2\pi^2}{(m+1)^3} + O(\frac{1}{m^4})}\\& = e^{-\frac{2\pi^2}{x_0^3 N} + O(N^{-1-\delta}) + O(N^{-\frac 43})} = e^{-\frac{\mu}{N} + O(N^{-1-\delta})}\,,\nonumber
		\end{align}
		which together with
		$$
		\frac{m+1}{m}\frac{\tan^2\frac{\pi}{m}}{\tan^2\frac{\pi}{m+1}} = 1+ O(N^{-\frac 13})
		$$
		and \eqref{eq:A2} yields
		\bb\label{eq:prod2}
		\prod_{t=1}^R\frac{A_{m-1,N_t}}{A_{m,N_t}} = e^{-\mu\frac{N_1+\dots+N_R}{N}}\big(1+O(N^{-\delta})\big) = e^{-\mu}\big(1+O(N^{-\delta})\big)\,.
		\ee
		Inserting \eqref{eq:prod1} and \eqref{eq:prod2} into \eqref{defVN}, the sum over $N_1,\dots, N_R$ can be performed and yields a combinatorial factor 
		\begin{align}\nonumber
			\binom{N-N'-R\lfloor \epsilon N\rfloor +R-1}{R-1} & \geq \frac{N^{R-1}}{(R-1)!}\Big(1-\frac{N'+R\lfloor\epsilon N\rfloor}{N}\Big)^{R-1} \\ &\geq  \frac{N^{R-1}}{(R-1)!}(1- \epsilon R)^{R-1}\big(1+O(N^{-1})\big)\,.\nonumber
		\end{align}
		Thus, we obtain
		\begin{align}\nonumber
			V_N \geq 4^{N-N'}\frac{1-e^{-\mu}}{(R-1)!}\sum_{|x_m-x_0|\leq N^{-\delta}}& N^{R-1} e^{-\mu m}\Big(\frac{1}{m+1}\tan^2\frac{\pi}{m+1}\Big)^R \\ & \cdot\big(1+ \tan^2\frac{\pi}{m+1}\big)^{-(N-N')} (1-\epsilon K)^{K} \big(1+O(N^{-\delta})\big)\,.\nonumber
		\end{align} 
		Using now \eqref{tanmscale}, the sum can be estimated by repeating the arguments leading to  \eqref{bound3}, 
		and we arrive at 
		\bb\nonumber
		V_N \geq 4^{N-N'}\frac{e^{\mu}-1}{(R-1)!}\sqrt{\frac{\pi}{B}}\Big(\frac{\pi^2}{x_0^3}\Big)^R e^{-A N^{\frac 13}}N^{-\frac 56} \cdot (1-\epsilon K)^{K} \big(1+O(N^{\frac 13-3\delta})\big)\,,
		\ee
		provided $\delta >\frac 19$. 
		Using this estimate in \eqref{bound0} as well as Theorem \ref{thm:4.1} and \eqref{x0mu}, the claimed inequality \eqref{est1} follows from \eqref{Rsum}.
	\end{proof}
	
	In the following, let $\Xi(T_0;M,\epsilon)$ denote the large-$N$ limit of the right-hand side of \eqref{est1},
	$$
	\Xi(T_0;M,\epsilon) =  \frac{e^{-\mu(r-1)}}{4^{|T_0|-K}} \sum_{R=1}^{K} \binom{K}{R}\frac{1}{(R-1)!}\Big(\frac{\mu}{2}\Big)^{R-1}\Big(\sum_{S=1}^M C_{S-1}4^{-S}\Big)^{K-R}(1-\epsilon K)^K\,,
	$$
	and set 
	\bb\label{defXi}
	\Xi(T_0) := \lim_{\substack{\epsilon\to 0\\ M\to \infty}} \Xi(T_0,M,\epsilon) =  \frac{e^{-\mu(r-1)}}{4^{|T_0|-K}} \sum_{R=1}^{K} \binom{K}{R}\frac{1}{(R-1)!}\Big(\frac{\mu}{2}\Big)^{R-1}2^{R-K}\,,
	\ee
	where we have also used that 
	$$
	\sum_{S=1}^\infty C_{S-1} 4^{-S} = X(\frac 14) = \frac 12\,.
	$$
	We then have following analogue of Lemma \ref{lem:addupto}.
	
	\begin{lem}\label{lem:sumrule2} For all $r\geq 1$ it holds that
		\bb\label{sumrule2}
		\sum_{T_0: h(T_0)=r} \Xi(T_0) = 1\,.
		\ee
	\end{lem}
	
	\begin{proof} 
		We use induction with respect to $r$. The case $r=1$ being trivial, let $r\geq 2$ and assume \eqref{sumrule2} holds for $r-1$. 
		It is convenient to rewrite \eqref{defXi} as 
		$$
		\Xi(T_0) = \sum_{\substack{B\subseteq D_r(T_0)\\ B\neq \emptyset}} \Xi(T_0,B)\,,
		$$
		where 
		$$
		\Xi(T_0,B) =  \frac{e^{-\mu(r-1)}}{4^{|T_0|-|D_r(T_0)|}}\frac{1}{(|B|-1)!}\Big(\frac{\mu}{2}\Big)^{|B|-1}2^{-|D_r(T_0)\setminus B|}\,,
		$$
		such that the sum in \eqref{sumrule2} becomes a double sum over $(T_0,B)$. 
		
		Let  $T'_0\in\cT_{\rm fin}$ be a fixed tree of height $r-1$ with $K'$ vertices at height $r-1$. Given a non-empty $B'\subseteq D_{r-1}(T_0')$, we then consider the contribution to the sum in \eqref{sumrule2} from all $(T_0,B)$ such that $T_0$ coincides with $T_0'$ up to height $r-1$ and such that $B'$ is the set of parents to vertices in $B$, i.e. $\phi_r(B)=B'$ in the notation of section \ref{sec:2}, where $\phi_r$ is the $r$'th parent map of $T_0$. We claim that this contribution is precisely  $\Xi(T'_0,B')$.
		
		In order to establish the claim, we first note that for any ordered set $B$ with $|B|=R$ the number of surjective order preserving maps $\varphi: B\to B'$ equals $ \binom{R-1}{ R'-1}$, where $R'=|B'|$. Moreover, for any given such $\varphi$, the number of ways of extending it to an order preserving map $\phi_r$ (not necessarily surjective) from any ordered set of $K$ elements, containing $B$ as an ordered subset, into $D_{r-1}(T_0')$ is easily seen to equal $\binom{ K+K'-1}{ K'+R-1}$. Hence, we have 
		$$
		\sum_{\substack{(T_0,B): h(T_0)=r,\, B\subseteq D_r(T_0)\\ B_{r-1}(T_0)=T_0',\, \phi_r(B)=B'}} \Xi (T_0,B) =  \frac{e^{-\mu(r-1)}}{4^{|T_0'|}} \sum_{K=1}^\infty\sum_{R=1}^K  \binom{K+K'-1}{K'+R-1} \binom{R-1}{R'-1}\frac{(\frac{\mu}{2})^{R-1}}{(R-1)!}2^{R-K}.
		$$
		Interchanging the summation order and using the identity \eqref{combid},
		the right-hand side becomes 
		\begin{align}\nonumber
			\frac{e^{-\mu(r-1)}}{4^{|T_0'|}} \sum_{R=R'}^\infty \binom{R-1}{R'-1}\frac{\mu^{R-1}}{(R-1)!}\,2^{K'+1}
			&=  \frac{e^{-\mu(r-1)}}{4^{|T_0'|}} \sum_{R=R'}^\infty \frac{\mu^{R-1}}{(R'-1)!(R-R')!}2^{K'+1}\\\nonumber
			&=  \frac{e^{-\mu(r-2)}}{4^{|T_0'|}} \frac{\mu^{R'-1}}{(R'-1)!}\,2^{K'+1}\,,
		\end{align}
		which is seen to be equal to $\Xi(T'_0,B')$, as claimed.
		
		Using this result, we get 
		\begin{align}\nonumber
			\sum_{\substack{(T_0,B): h(T_0)=r\\ B\subseteq D_r(T_0),\, B\neq \emptyset}} \Xi(T_0,B) &= \sum_{\substack{(T'_0,B'): h(T'_0)=r-1\\ B'\subseteq D_r(T'_0),\, B'\neq \emptyset}} \sum_{\substack{(T_0,B): h(T_0)=r, B\subseteq D_r(T_0) \\ B_{r-1}(T_0)=T'_0, \phi_r(B)=B'}} \Xi(T_0,B)\\& =  \sum_{\substack{(T'_0,B'): h(T'_0)=r-1\\ B'\subseteq D_r(T'_0),\, B'\neq \emptyset}} \Xi(T'_0,B') =1\,,\nonumber
		\end{align}
		where the induction assumption has been used in the last step. This finishes the proof.
	\end{proof}
	
	On the basis of Lemma \ref{lem:sumrule2},  the arguments showing tightness in the proof of  Corollary \ref{tightness0} can be repeated, thus establishing the following corollary. Details of the proof are left to the reader.
	
	\begin{cor}\label{cor:tightness2}
		For each $\mu>0$, the sequence of measures $(\nu^{(\mu)}_N)$ on $\cT$ is tight\,. 
	\end{cor}
	
	Similarly, we obtain the following main result of the existence of the local limit.
	
	\begin{thm}\label{thm:limit2}
		For each $\mu>0$, the sequence $(\nu^{(\mu)}_N)$ is weakly convergent to a Borel probability measure $\nu^{(\mu)}$ on $\cT$ characterized by
		\bb\label{ballmeasure2}
		\nu^{(\mu)}(\cB_{\frac 1r}(T_0)) = \Xi(T_0) = \frac{e^{-\mu(r-1)}}{4^{|T_0|}}2^{K+1}\ \sum_{R=1}^{K} \binom{K}{R}\frac{\mu^{R-1}}{(R-1)!}\,,
		\ee 
		for any tree $T_0\in\cT_{\rm fin}$ of height $r\geq 1$, where $K=|D_r(T_0)|$.
	\end{thm}
	
	\begin{proof} This follows by the same line of reasoning as in the proof of Theorem \ref{thm:limit1}, using  Corollary \ref{cor:tightness2}, Lemma \ref{lem:sumrule2}, and  Lemma \ref{lem:4.2}. We leave the details to the reader. 
	\end{proof}

	\subsection{Properties of the local limit}\label{sec:4.3}
	
	Before formulating the basic decomposition result for the measure $\nu^{(\mu)}$, some further notational conventions are needed.  
	
	Given $n\in\mathbb N$, let the standard $(n-1)$-simplex $\Delta_{n}$ be defined by
	$$
	\Delta_{n} := \{ (x_1,\dots,x_n) \mid x_1 + \dots +x_n = 1\,,\; x_1,\dots, x_n > 0 \}\,.
	$$
	Scaling by $\mu>0$ we obtain the simplex
	$$
	\mu\cdot \Delta_{n} := \{ (\mu_1,\dots,\mu_n) \mid \mu_1 + \dots +\mu_n = \mu\,,\; \mu_1,\dots,\mu_n > 0 \}\,.
	$$
	By $d\omega(\mu_1,\dots,\mu_n)$ we shall denote the normalised Lebesgue measure on $\mu\cdot\Delta_n$.
	
	For $(\mu_1,\dots,\mu_n) \in \mu\cdot \Delta_{n}$, consider the product measure $\prod_{k=1}^n \nu^{(\mu_k)}$ on $\cT^n$ and note that, if $\cB_1,\dots,\cB_n$ are balls in $\cT$, then $\prod_{k=1}^n \nu^{(\mu_k)}(\cB_1\times\dots\times\cB_n)$ is a continuous function of $(\mu_1,\dots,\mu_n)$ by \eqref{ballmeasure2}.  For $r\in\mathbb N$, let  ${\cal F}_r^n,\, r\in\mathbb N$,  denote the collection of subsets of $\cT^n$ that can be written as countable unions of products of balls in $\cT$ with radius $\frac 1r$ and set ${\cal F}_\infty^n = \cup_{r\in\mathbb N} {\cal F}_r^n$. As in the case $n=1$ discussed in Remark \ref{remcT} iii), $\cF_\infty^n$ is a set algebra consisting of sets with empty boundary, and it generates the Borel $\sigma$-algebra $\cF^n$ of $\cT^n$. It follows that $\prod_{k=1}^n \nu^{(\mu_k)}(A)$ is a measurable function of $(\mu_1,\dots,\mu_n)$ for any $A\in\cF_\infty^n$, and we obtain a well defined finitely additive set function $\lambda_0$ on $\cF_\infty^n$
	by setting  
	\bb\label{defintmeas}
	\lambda_0 (A) := \int_{\mu\cdot \Delta_{n}}\prod_{k=1}^n \nu^{(\mu_k)}(A)\,d\omega(\mu_1,\dots,\mu_n)\,,\quad A\in {\cal F}_\infty^n\,.
	\ee
	Moreover, it follows from the monotone convergence theorem that $\lambda_0$ is countably additive on $\cF_\infty^n$, and hence extends to a unique Borel probability measure  on $\cT^n$ (see e.g. Theorem A in \S 13 in \cite{halmos1950measure}). We shall denote this extension by 
	$$
	\int_{\mu\cdot\Delta_{n}}d\omega \prod_{k=1}^n \nu^{(\mu_k)}\,.
	$$
	The right-hand side of eq. \eqref{defintmeas} then equals the measure of $A$ for more general sets, such as open and closed sets, but we shall only need this expression for products of balls.
	
	In the following theorem we make use of the measures just introduced as well as of the BGW measure $\rho$ introduced in Remark \ref{muequalzero}. Note that, under the identification of $ \cB_{\frac 1r}(T_0)$ with $\cT^{D_r(T_0)|}$ implied by Lemma \ref{lem:2.3} ii) where $T_0\in\cT_{\rm fin}$ has height $r$, the identity 
	\bb\label{rhocond}
	\rho(\cdot\mid  \cB_{\frac 1r}(T_0)) = \prod_{i\in D_r(T_0)} \rho_i\,,\quad \rho_i=\rho\,,
	\ee
	holds  as a consequence of \eqref{rhoball} and \eqref{def:rho}, expressing that in the corresponding BGW process the events of birth by different individuals in any generation are i.i.d. according to $\rho$.
	
	\begin{thm} \label{decomposition}
		Assume $\mu>0$ and that $T_0 \in \cT_{\rm fin}$ has height $r$, and let $K=|D_r(T_0)|$. Under the identification $\mathcal{B}_{\frac{1}{r}} (T_0) = \cT^{D_r(T_0)}$  provided by Lemma \ref{lem:2.3} ii), the measure $\nu^{(\mu)}$ restricted to $\mathcal{B}_{\frac{1}{r}} (T_0)$ is given by
		\begin{align} \nonumber
			\nu^{(\mu)} \big|_{\mathcal{B}_{\frac{1}{r}}(T_0)} &= e^{-(r-1)\mu} 4^{-|T_0|}2^{K+1} \sum\limits_{\substack{D \subseteq D_r(T_0)\\ D \neq \emptyset}} \frac{\mu^{|D|-1}}{(|D|-1)!} \\
			& \Big(\int\limits_{\mu\cdot\Delta_{|D|}} d\omega \prod\limits_{i \in D} \nu^{(\mu_i)}\Big) \times  \prod\limits_{j \in D_r(T_0) \setminus D} \rho_j \label{decomp}\,,
		\end{align}
		where $\rho_j = \rho$.
	\end{thm}
	
	\begin{proof}
		It should be possible to give a proof based on the result of Theorem \ref{thm:limit2} alone. Here we give a limit argument using estimates similar to those applied in the proofs of Theorem \ref{thm:4.1} and Lemma \ref{lem:4.2}. It is sufficient to show that both sides of equation \eqref{decomp} coincide when acting on sets of the form $\mathcal{B}_{\frac{1}{h_1}} (S_1)\times \dots \times \mathcal{B}_{\frac{1}{h_K}}(S_K)$ for arbitrary finite trees $S_1, \dots , S_K$ where $h_i = h(S_i)$. Thus, let $S_1,\dots,S_K$ be given and let $T_0^\prime$ denote the tree obtained by grafting them onto the vertices in $D_r(T_0)$ in the given order from left to right. Moreover, let
		$$
		l_i = |D_{h_i} (S_i)| \quad\mbox{and}\quad 		L= l_1 + \dots + l_K\, .
		$$
		By the definition of the grafting process, any tree corresponding to an element in $\mathcal{B}_{\frac{1}{h_1}} (S_1)\times \dots \times \mathcal{B}_{\frac{1}{h_K}}(S_K)$ is then obtained by grafting  trees $T_1, \dots , T_L \in \cT$ onto the vertices (leaves) of $T_0^\prime$ that belong to $\cup_{i=1}^K D_{h_i}(S_i)$. Calling the set of trees obtained in this way  $\mathcal{B}(T_0 ^\prime)$, we hence need to compute its  $\nu ^{(\mu)} _N$-measure and evaluate its large-$N$ limit. For this purpose we use a procedure similar to the one applied in the proof of Lemma \ref{lem:4.2}, taking into account that the vertices onto which the trees $T_1,\dots, T_L$ are grafted are not at equal height. 
		
		Let us write $T_1,\dots,T_L$ as $T_1 ^{(1)},\dots, T^{(1)} _{l_1} ,\dots, T_1 ^{(K)}, \dots, T^{(K)} _{l_K}$, such that the trees $T_1^i,\dots,T_{l_i}^i$ are grafted from left to right onto the vertices in $D_{h_i}(S_i)$, and set 
		$$
		N(i):= \sum\limits_{j=1} ^{l_i} |T^{(i)} _j|\,.
		$$
		As in the proof of Lemma \ref{lem:4.2}, we  call $T^{(i)} _{j}$ \textit{small} if $|T^{(i)} _{j}| \leq M$, while we call it $large$ if 
		\bb\label{large2}
		|T^{(i)} _j| > \epsilon N(i)\quad \mbox{and} \quad N(i) > \epsilon N\,.
		\ee
		Assuming that $N$ fulfills the inequalities $N > \frac{M}{\epsilon^2}$ and $N>L \cdot M + |T_0 ^\prime| -L$ it is seen  that no $T^{(i)} _j$ can be both small and large and there must be at least one large branch $T^{(i)} _j$. The set of vertices onto which the large branches among $T_1 ^{(i)} ,\dots, T_{l_i}^{(i)}$ are grafted will be denoted by $B_i \subseteq D_{h_i} (S_i)$, and we set
		$$
		N^{(i)} := \sum\limits_{v \in B_i} N_v ^{(i)}\,,
		$$
		where $N_v ^{(i)}$ is  the size of the tree grafted onto a vertex $v \in D_{h_i} (S_i)$. Moreover,  the total size of small branches is set to 
		$$ 
		N^\prime :=\sum\limits_{(i,v): v \notin B_i} N_v ^{(i)}\,.
		$$
		
		We then have
		
		\begin{align}
			\nonumber \nu^{(\mu)} _N (\mathcal{B}(T_0 ^\prime)) & \geq \sum\limits_{\substack{B_i \subseteq D_{h_i} (S_i) \\ i = 1,\dots,K}} \frac{e^{-(r-1)\mu}}{Z^{(\mu)}_N} \sum\limits_{|x_m-x_0|\leq N^{-\delta}} e^{-m\mu} \sum \limits_{N_v^{(i)} \leq M , v \notin B_i} \prod\limits_{i=1} ^K \prod\limits_{v \notin B_i} C_{N_v^{(i)} -1} \\
			& \sum\limits_{\substack{N(i)> \epsilon N,~ N_v ^{(i)} > \epsilon N(i), ~v \in B_i \\ \sum\limits_{i} N^{(i)} = N- N^\prime +L - |T_0^\prime|}} \Big( \prod\limits_{i=1} ^K \prod\limits_{v \in B_i} A_{m-h_i+1,N_v ^{(i)}} - \prod\limits_{i=1} ^K \prod\limits_{v \in B_i}A_{m-h_i,N_v^{(i)}} \Big) \,, \label{dert1}
		\end{align}
		where $x_m$ and $x_0$ are defined by \eqref{eq:xm} and \eqref{eq:x0} and $\delta\in ]0,\frac 16[$.
		Rewriting the last parenthesis as
		$$
		\Big( \prod\limits_{i=1}^K \prod\limits_{v\in B_i}  A_{m, N_v ^{(i)}} \Big)\Big( \prod\limits_{i=1}^K \prod\limits_{v\in B_i} \frac{A_{m-h_i+1,N_v^{(i)}}}{A_{m,N_v^{(i)}}} \Big) \Big( 1- \prod\limits_{i=1}^K \prod\limits_{v\in B_i} \frac{A_{m-h_i,N_v^{(i)}}}{A_{m-h_i+1,N_v ^{(i)}}} \Big)\,,
		$$
		we have by a similar computation as the one leading to \eqref{eq:prod2} that the last parenthesis in this expression equals
		\bb \label{1A}
		(1-e^{-\mu} )\big(1+ O(N^{-\delta})\big)\,,
		\ee
		\noindent while the second one equals
		\bb \label{2A}
		\prod\limits_{i=1} ^K e^{-\mu(h_i-1)\frac{N^{(i)}}{N}} \big(1+ O(N^{-\delta})\big)\,.
		\ee
		Applying the saddle point method in the same manner as in the proof of Theorem \ref{thm:4.1} we get, using identical notation,
		\begin{align} \label{3A}
			\sum\limits_{|x_m-x_0|\leq N^{-\delta}} e^{-\mu m}  \prod\limits_{i=1}^K \prod\limits_{v\in B_i}  A_{m, N_v ^{(i)}} \geq \frac{4^{\sum\limits_{i}N^{(i)}}}{N^{\sum\limits_i |B_i|}}\, e^\mu \prod\limits_{i=1} ^K \Big(\frac{\mu}{2}\Big)^{|B_i|}\sqrt{\frac{\pi}{B}}\, N^{\frac{1}{6}} \,e^{-AN^{\frac{1}{3}}}\big(1+O(N^{\frac{1}{3}-3\delta})\big) \, .
		\end{align}
		Using \eqref{1A},  \eqref{2A} and \eqref{3A} in  \eqref{dert1}, the summation over the sizes $N_v ^{(i)}, v\in B_i$, of the large branches for fixed $N(i)$  
		yields, for each $i=1,\dots, K$,  a combinatorial factor
		$$
		\binom{N^{(i)} - \lfloor\epsilon N(i)\rfloor |B_i|+ |B_i|-1}{|B_i|-1} \geq  \frac{(N^{(i)})^{|B_i|-1}}{(|B_i|-1)!} (1-\epsilon L)^{|B_i|-1}\big(1+O(N^{-1})\big) \,,
		$$
		where the inequality follows by observing that $|B_i|\leq L$ and that 
		$$
		N(i) \geq N^{(i)} \geq N(i) - LM \geq N(i)\big(1- \frac{LM}{\epsilon N}\big)\,,\quad\mbox{if $B_i\neq\emptyset$}\,,
		$$
		as a consequence of condition \eqref{large2}.
		
		Setting
		$$
		D=\{ i  \big| ~ B_i \neq \emptyset \}\,,
		$$
		considered as a subset of $D_r(T_0)$, we thus obtain
		\begin{align}
			\nonumber & \nu^{(\mu)}_N(\mathcal{B}(T_0 ^\prime)) \geq  \frac{e^{-(r-1)\mu}}{Z^{(\mu)}_N} (e^\mu -1) \sqrt{\frac{\pi}{B}} \, N^{\frac{1}{6}}\, e^{-AN^{\frac{1}{3}}} \sum\limits_{\substack{B_i \subseteq D_{h_i} (S_i)\\i=1,\dots,K}} \Big[ \sum\limits_{N_v ^{(i)} \leq M, v \notin B_i} \prod\limits_{(i,v): v \notin B_i} C_{N_v ^{(i)} -1} \Big]  \\
			& \sum\limits_{\substack{N^{(i)}> \epsilon N,\, i\in D \\ \sum\limits_{i\in D} N^{(i)}=N-N^\prime +L-|T_0 ^\prime|}} N^{-|D|}\prod_{i\in D} 4^{N^{(i)}}  e^{-\mu(h_i-1)\frac{N^{(i)}}{N}} \Big(\frac{\mu}{2} \Big)^{|B_i|} \frac{\Big((1-\epsilon L)\frac{N^{(i)}}{N}\Big)^{|B_i|-1}}{(|B_i|-1)!} \nonumber\\ &~~~~~~~~~~~~~~~~~~~~~~~~~~~~~~~~~~~~~~~~~~~~~~~~~~~~~~~~~~~~~~~~~~~~~~~~~~~~\cdot~\big(1+ O(N^{\frac{1}{3}-3\delta})\big)\,. \label{dert2}
		\end{align}
		
		Observing that the constraint on the summation variables $N^{(i)}$ implies that
		$$
		\sum_{i\in D} \frac{N^{(i)}}{N} = 1 + O(N^{-1})\,,
		$$
		it is seen that, up to a factor $N$, the last sum in \eqref{dert2} can be viewed as a Riemann sum approximating an integral over $\Delta_{|D|}$. More precisely, since the integrand below is a smooth function on $\mathbb R^{|D|}$, we have  
		\begin{align*}
			&\sum\limits_{\substack{N^{(i)}\geq \epsilon N,\, i\in D \\ \sum\limits_{i \in D} N^{(i)} = N-N^\prime +L- |T_0^\prime|}}  \prod\limits_{i\in D}  e^{-\mu (h_i-1) \frac{N^{(i)}}{N}} \Big( \frac{N^{(i)}}{N} \Big)^{|B_i|-1} N^{-(|D|-1)} \\
			&=\frac{1}{(|D|-1)!} \int\limits_{\Delta_{|D|} ^{(\epsilon)}} \prod \limits_{i \in D} e^{-\mu(h_i-1)x_i} x_i^{|B_i|-1}d\omega(x_1,\dots,x_{|D|})(1+ O(N^{-1}))\,,
		\end{align*}
		where
		$$
		\Delta_{|D|} ^{(\epsilon)} = \{ x=(x_1,\dots,x_{|D|}) \big| ~x_1+\dots+x_{|D|} =1, ~x_i > \epsilon \}\,.
		$$
		Inserting this into \eqref{dert2} and taking into account \eqref{eq:asympZ}, we hence get
		\begin{align}
			\nonumber &\nu^{(\mu)}_N (\mathcal{B}(T_0 ^\prime)) \geq \sum\limits_{\substack{D \subseteq D_r(T_0)\\ D \neq \emptyset}} \sum\limits_{\substack{B_i \subseteq D_{h_i} (S_i)\\i \in D,B_i \neq \emptyset}} e^{-(r-1)\mu} 4^{-|T_0 ^\prime| +L} \big( \frac{\mu}{2} \big)^{-1} \Big( \sum_{S=1}^M C_{S-1}4^{-S}\Big) ^{L-\sum_{i\in D}|B_i|} \\\label{dert3}
			& \cdot\frac{1}{(|D|-1)!} \int\limits_{\Delta^{(\epsilon)}_{|D|}} d\omega(x) \prod\limits_{i \in D} \big( \frac{\mu}{2} \big)^{|B_i|} e^{-\mu(h_i-1)x_i} \frac{\big((1-\epsilon L)x_i\big)^{|B_i|-1}}{(|B_i|-1)!}\big(1+O(N^{\frac 13 -3\delta})\big)\,.
		\end{align}
		Since $\mathcal{B}(T_0 ^\prime)$ has empty boundary by Lemma \eqref{lem:2.3} i), the limit of the left-hand side of \eqref{dert3} as $N\to\infty$ equals $\nu^{(\mu)}(\mathcal{B}(T_0 ^\prime))$. Thus, taking first the large-$N$ limit and subsequently letting $M\to\infty$ and $\epsilon\to 0$,  we arrive at the inequality
		\begin{align}
			\nonumber &\nu^{(\mu)} (\mathcal{B}(T_0 ^\prime)) \geq \sum\limits_{\substack{D \subseteq D_r(T_0)\\ D \neq \emptyset}} \sum\limits_{\substack{B_i \subseteq D_{h_i} (S_i)\\i \in D,B_i \neq \emptyset}} e^{-(r-1)\mu} 4^{-|T_0 ^\prime| +L} \big( \frac{\mu}{2} \big)^{-1} \Big( \frac 12\Big) ^{L-\sum_{i\in D}|B_i|} \\\nonumber
			& \cdot\frac{1}{(|D|-1)!} \int\limits_{\Delta_{|D|}} d\omega(x) \prod\limits_{i \in D} \big( \frac{\mu}{2} \big)^{|B_i|} e^{-\mu(h_i-1)x_i} \frac{x_i^{|B_i|-1}}{(|B_i|-1)!}\,,
		\end{align}
		which by a scaling of integration variables and setting $|B_i|=R_i$ can be rewritten as 
		\begin{align}
			\nonumber
			\nu^{(\mu)} (\mathcal{B}(T_0 ^\prime))&\geq e^{-(r-1)\mu} 4^{-|T_0|} 2^{K+1} \sum\limits_{\substack{D \subseteq D_r(T_0)\\ D \neq \emptyset}}\frac{\mu^{|D|-1}}{(|D|-1)!} \prod\limits_{i \notin D} 4^{-|S_i|} 2^{l_i+1} \\
			\int\limits_{\mu\Delta_{|D|}} &d\omega(\mu_1,\dots,\mu_{|D|}) \prod\limits_{i \in D} \sum\limits_{R_i=1} ^{l_i} \binom{l_i}{R_i}e^{-(h_i-1)\mu_i}\frac{\mu_i^{R_i-1}}{(R_i-1)!}\, 4^{-|S_i|}\, 2^{l_i+1}\,. \label{dert4}
		\end{align}
		Recalling \eqref{ballmeasure2} and the definition of the BGW measure $\rho$, it follows that the last expression precisely equals the right-hand side of \eqref{decomp}, when applied to $\mathcal{B}_{\frac{1}{h_1}} (S_1)\times \dots \times \mathcal{B}_{\frac{1}{h_K}}(S_K)$. Since the integral on the right-hand side of \eqref{decomp} is a probability measure for each nonempty $D\subseteq D_r(T_0)$, the total mass of the measure on the right-hand side equals 
		$$
		e^{-(r-1)\mu} 4^{-|T_0|} 2^{K+1} \sum\limits_{\substack{D \subseteq D_r(T_0)\\ D \neq \emptyset}}\frac{\mu^{|D|-1}}{(|D|-1)!}\,,
		$$
		which by \eqref{ballmeasure2} equals $\nu^{(\mu)}(\cB_{\frac 1r}(T_0))$. This shows that summing the two sides of \eqref{dert4} over trees $S_1,\dots, S_K$ of heights $h_1,\dots, h_K$, respectively, yields identical results. Therefore we conclude that equality must hold in \eqref{dert4}, which completes the proof of \eqref{decomp}.
	\end{proof}

	\begin{remark} Using the identity
		$$
		\int_{\Delta_n} d\omega(x) \prod_{i=1}^nx_i^{m_i}  = \frac{(n-1)!}{(m_1+\dots+m_n +n-1)!}\prod_{i=1}^n m_i!\,,
		$$
		for non-negative integers $m_1,\dots, m_n$,
		the reader may easily verify that if $h_1=\dots=h_n :=h$, in which case $\cB(T_0^\prime)=\cB_{\frac{1}{r+h-1}}(T_0^\prime)$, formula \eqref{decomp} (or the right-hand side of \eqref{dert4}) reproduces the expression \eqref{ballmeasure2} in this case. 
	\end{remark}
	
\begin{figure}[H]
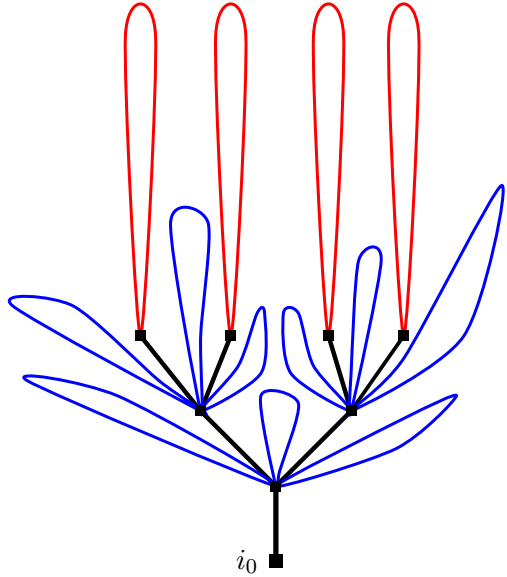

		\centering
		\begin{align*}
			\raisebox{0cm}{\treespine} 
		\end{align*}
		\captionof{figure}{Structure of the elements in $\chi^{-1}(\cB_{\frac 13} ^s (T_0 ^s))$ where $T_0 ^s$  of height $h(T_0 ^s)=3$ with root vertex $i_0$ is shown in black. Red blobs and blue blobs indicate  the infinite branches and the critical BGW branches, respectively.}
		\label{figure:treespine}
	\end{figure}

	As a first consequence of Theorem \ref{decomposition} we note that, since $\nu^{(\mu)}$ is supported on $\cT_\infty$ while $\rho$ is supported on $\cT_{\rm fin}$, formula \eqref{decomp}  provides an explicit decomposition of $\nu^{(\mu)}$ conditioned on the event $\{T\mid B_r(T) = T_0\}$ into measures supported on trees that have a fixed spine up to height $r$ labelled by the set $D$ of vertices of infinite type at height $r$. In order to elaborate further on this decomposition, we first determine the measure $\tilde\nu^{(\mu)}$ induced by $\nu^{(\mu)}$ on the set $\cT^s$ of spine trees, i.e. $\tilde\nu^{(\mu)}$  is the pushforward of $\nu^{(\mu)}$ by the spine map $\chi$, 
	\bb\label{deftildenu}
	\tilde\nu^{(\mu)}(A) = \nu^{(\mu)}\big(\chi^{-1}(A)\big)\,,\quad A\subseteq\cT^s\, \mbox{a Borel set}\,.
	\ee
	We shall use the notation 
	$$
	\cB^s_{a}(T^s) := \cB_a(T^s)\cap \cT^s
	$$ 
	for the ball in $\cT^s$ of radius $a>0$ around a spine tree $T^s$. Given $T^s$ and $r\in\mathbb N$,  we have by Remark \ref{remcT} i) that $\cB^s_{\frac 1r}(T^s) = \cB_{\frac 1r}(T^s_0)\cap \cT^s$, where $T_0^s =B_r(T^s)$ is a finite tree of height $r$, all of whose leaves are at height $r$.   Moreover, a tree $T$ in $\chi^{-1}(\cB^s_{\frac{1}{r}}(T^s))$ is obtained by grafting arbitrary  infinite trees onto the leaves of $T^s_0$ while grafting finite trees $T_{(i,n)}$ in each sector $(i,n)$ of any vertex at height $<r$ different from the root, see Fig.\ref{figure:treespine}. In this way, we have by Lemma \ref{lem:2.3} i) a homeomorphism
	\bb\label{id}
	G: \chi^{-1}(\cB^s_{\frac{1}{r}}(T^s)) \to \cT_\infty^R\times  \cT_{\rm fin}^\sigma\,,
	\ee
	where $R=|D_r(T^s)|$ is the number of leaves in $T_0^s$ and 
	\bb\label{sigmaform}
	\sigma = 2|T^s_0| -R-1
	\ee
	is the total number of sectors associated with $T^s$, excluding those adjacent to the leaves.  
	
	\begin{thm} \label{backbonedist}
		Assume $\mu>0$ and let $T^s_0$ be any finite tree of height $r$ with $R$ leaves, all of which are at height $r$. Then the following two statements hold.
		
		i)\; 	\bb \label{nutilde}
		\tilde{\nu} ^{(\mu)} (\cB^s_{\frac{1}{r}}(T^s_0))=e^{-(r-1)\mu} \frac{\mu^{R-1}}{(R-1)!} ~,~~ r \geq 1\,.
		\ee
		
		ii)\; Under the identification \eqref{id} we have 
		\bb\label{decomp2}
		\nu^{(\mu)}(\cdot\mid \chi^{-1}(\cB^s_{\frac{1}{r}}(T_0^s))) = \Big(\int_{\mu\Delta_R}d\omega\prod_{j=1}^R\nu^{(\mu_i)}\Big)\times \prod_{(i,n)} \rho_{(i,n)}\,,
		\ee
		where $\rho_{(i,n)}=\rho$. Equivalently, conditioning on the event $\{T\mid B_r(\chi(T)) = T_0^s\}$ renders the finite branches $T_{(i,n)}$ i.i.d. according to $\rho$, whereas the joint distribution of the infinite branches is independent thereof and determined by $\int_{\mu\Delta_R}d\omega\prod_{j=1}^R\nu^{(\mu_i)}$. 
	\end{thm}
	
	\begin{proof} In order to compute the left-hand side of \eqref{nutilde}, we use the preceding notation and set $T^\prime_{(i,n)}=B_{r-h_i+1}(T_{(i,n)})$, where $h_i$ is the height of $i$ in $T^s_0$. Letting $T_0$ be the tree of height $r$ obtained by grafting $T^\prime_{(i,n)}$ onto $T^s_0$ at $(i,n)$ instead of $T_{i,n}$ in each sector, Theorem \ref{decomposition} shows, for fixed $T^\prime_{(i,n)}$, that the branches of $T_{(i,n)}$ grafted on the vertices of $T^\prime_{(i,n)}$ at (maximal) height $r-h_i+1$ are i.i.d. according to $\rho$. If $l_{(i,n)}$ denotes the number of vertices in $T^\prime_{(i,n)}$ at height $r-h_i+1$, the number of vertices in $T_0$ at height $r$ equals $K=R+\sum_{(i,n)} l_{(i,n)}$. Together with \eqref{sigmaform} this implies, with notation as in Theorem \ref{decomposition}, that 
		$$
		4^{-|T_0|}\,2^{K+1}  =  4^{-|T^s_0|}\,2^{R+1}\prod_{(i,n)} 4^{-|T^\prime_{i,n}|+1}\,2^{l_{(i,n)}} = \prod_{(i,n)} 4^{-|T^\prime_{i,n}|}\,2^{l_{(i,n)}+1}\,.
		$$
		Taking into account \eqref{rhoball} and \eqref{rhocond} and  the identification \eqref{id}, it then follows from Theorem \ref{decomposition} that
		$$
		\nu^{(\mu)}\Big|_{\chi^{-1}(\cB^s_{\frac{1}{r}}(T_0^s))} = e^{-(r-1)\mu} \frac{\mu^{R-1}}{(R-1)!}\Big(\int_{\mu\Delta_R}d\omega\prod_{j=1}^R\nu^{(\mu_i)}\Big)\times \prod_{(i,n)} \rho_{(i,n)}\,,
		$$
		which clearly implies both statements of the theorem.
	\end{proof}
	
	It is worth observing that the right-hand side of \eqref{nutilde} only depends on the height $r$ and the number of vertices at height $r$ in $T_0^s$, but not on the structure of $T_0^s$ below height $r$. As  a consequence it follows, in particular, that the random variables $|D_r(T^s)|,\, r\in\mathbb N$, constitute a Markov process. The following corollary shows that it is a random walk on the positive integers with Poisson distributed increments.
	
	\begin{cor} \label{cor:nutilde} The random variables  $\zeta_r,\,r\in \mathbb N$, on $\cT^s$ defined by
		$$
		\zeta_r(T^s) = |D_{r+1}(T^s)|-|D_{r}(T^s)|
		$$
		are i.i.d. according to  
		$$
		\tilde\nu^{(\mu)}(\zeta_r =n) =  e^{-\mu}\cdot \frac{\mu^{n}}{n!}\,,\quad n= 0,1,2,\dots.
		$$
		In particular,	the number $|D_r(T^s)|$ of vertices at height $r\geq 2$ with respect to $\tilde\nu^{(\mu)}$ is Poisson distributed with mean $(r-1)\mu$, i.e.
		\bb \label{poisson}
		\tilde\nu^{(\mu)} (\{|D_r(T^s)=R \} ) = e^{-(r-1)\mu}\frac{((r-1)\mu)^{R-1}}{(R-1)!}\,,\quad R\geq 1\, .
		\ee
	\end{cor}
	\begin{proof}
		By \eqref{nutilde} we have
		$$
		\tilde\nu^{(\mu)} (\{|D_r(T^s)|=R \} ) = e^{-(r-1)\mu} \frac{\mu^{R-1}}{(R-1)!} \cdot M_{R,r}\,,
		$$
		where $M_{R,r}$ is the number of finite trees $\tilde{T}_0$ with height $r$ and with $R$ leaves, all of which are at height $r$. This number can be determined recursively by noting that, 
		given $R \geq 1$ vertices at height $r$ and $R^\prime\geq 1$ vertices at height $r-1$, the number of ways of connecting them with non-intersecting edges such that every vertex at height $r-1$ has at least one offspring is given by $\binom{R-1}{R^\prime -1}$, interpreted as $0$ when $R<R^\prime$. Hence, 
		\eqref{nutilde} implies that
		\begin{align}\nonumber
			\tilde\nu^{(\mu)} (\{|D_r(T^s)|=R ,\; |D_{r-1}(T^s)|=R^\prime\} )& = e^{-(r-1)\mu} \frac{\mu^{R-1}}{(R^\prime -1)!(R-R^\prime)!} \cdot M_{R^\prime, r-1}\\\nonumber
			&= e^{-\mu}\frac{\mu^{R-R^\prime}}{(R-R^\prime)!} \tilde\nu^{(\mu)} (\{|D_{r-1}(T^s)|=R^\prime\})\,,
		\end{align}
		and consequently
		$$
		\tilde\nu^{(\mu)} (\{|D_r(T^s)|=R\} \mid \{|D_{r-1}(T^s)|=R^\prime\} ) 
		= e^{-\mu}\frac{\mu^{R-R^\prime}}{(R-R^\prime)!}\,.  
		$$ 
		Since $|D_{r}(T^s)|, \, r\in\mathbb N$, is a Markov process as already mentioned, the first statement of the corollary follows. The second one then follows by taking an $(r-1)$-fold convolution of the Poisson distribution \eqref{poisson} and using that $|D_1(T^s)|=1$. This completes the proof.
	\end{proof}
	
	We note that the measures $\nu^{(\mu)}$ also appear in \cite{abraham2020veryfat} as limits of conditioned BGW trees and some of its properties are discussed therein as well. 
	
	Denoting expectation w.r.t. $\tilde\nu^{(\mu)}$ by $\tilde{\mathbb E}_\mu$, the following result on the average volume growth of $T^s$ now easily follows.
	
	\begin{cor}\label{expectedbackbone}
		The following relations hold.
		\begin{itemize} 
			\item[i)]$\tilde{\mathbb E}_{\mu}(|D_r|)=\mu(r-1)$ ,
			\item[ii)] $\tilde{\mathbb E}_{\mu}(|B_r|)=\half \mu r(r-1) +1$,
			\item[iii)]  $\lim_{r\to\infty}\frac{|D_r(T^s)|}{r} = \mu $\quad\mbox{for $\tilde\nu^{(\mu)}$-a.e. $T^s$},
			\item[iv)]  $\lim_{r\to\infty}\frac{|B_r(T^s)|}{r^2} = \frac 12 \mu $\quad\mbox{for $\tilde\nu^{(\mu)}$-a.e. $T^s$}.
		\end{itemize}
	\end{cor}

	\begin{proof}
		The relation i)  follows immediately from \eqref{poisson}, and iii) follows by applying the strong law of large numbers (see e.g. \cite{breiman1968probability}). The second relation follows from i) and \eqref{DB}, and similarly iv) follows from iii).
	\end{proof}
	
	\begin{remark}
		Corollary \ref{expectedbackbone} iii) shows that, for $\mu>0$, it holds that $d_h=2$ for $\tilde{\nu}^{(\mu)}$-a.e. $T^s$. 
	\end{remark}
	
	Next we note the following result on the average volume growth of $T$ with respect to $\nu^{(\mu)}$.
	\begin{cor} \label{expected}
		The following statements hold.
		\begin{itemize}
			\item [i)] $\mathbb E_{\mu}(|D_r|)= \mu r^2\big(1+ O(\frac{1}{r})\big),$
			\item[ii)] $\mathbb E_{\mu}( |B_r|) = \frac{1}{3}\mu r^3\big(1+O(\frac{1}{r})\big)$.
		\end{itemize}
	\end{cor}
	\begin{proof}  It is clear that $ii)$ follows from $i)$ and \eqref{DB}. To establish i), we use that $\mathbb E_\rho(|D_r|) =1$ for all $r$, since $\rho$ is associated with a critical BGW process. Using this in \eqref{decomp2} together with 
		$$
		D_r(T) = D_r(\chi(T)) + \sum_{(i,n)} |D_{r-h_i+1}(T_{(i,n)})|\,,
		$$
		with notation as above,  we get
		$$
		{\mathbb E}_\mu (|D_r|\,\mid \chi^{-1}(\cB^s_{\frac{1}{r}}(T^s))) = \sigma(B_r(T^s) )+ |D_r(T^s)|\,, 
		$$
		where $\sigma$ is given by \eqref{sigmaform}. Integrating over $T^s$ then gives 
		$$
		\mathbb E_{\mu}(|D_r|) = \tilde{\mathbb E}_{\mu}(2 |B_r| -1)\,,
		$$
		and so relation i) follows from Corollary \eqref{expectedbackbone} ii). 
	\end{proof}
	
	We conclude by stating the following almost sure result on the volume growth of $T$ w.r.t. $\nu^{(\mu)}$, whose proof is deferred to the appendix.  
	
	\begin{thm}\label{asvolgrowth} For each $\mu>0$, there exist constants $C''_1, C''_2>0$ and for $\nu^{(\mu)}$-almost every $T\in\cT$ a number $r_0(T)\in\mathbb N$, such that 
		\bb\label{ashaus}
		C''_1\cdot r^3\;\leq\; |B_r(T)|\;\leq C''_2\cdot r^3 \log r\,, \quad r\geq r_0(T)\,.
		\ee
		In particular, it holds that $d_h=3$ for $\nu^{(\mu)}$-a.e. tree $T$.
	\end{thm}

	\subsection*{Acknowledgements}
	
	We thank Mireille Bousquet-Mélou and Nicolas Curien for bringing our attention to the papers \cite{guttmann2015analysis,abraham2020veryfat}. 
	
	\section*{Appendix}
	
	\subsection*{Proof of Proposition \ref{growth2}}	
	
	Note first, that the lower bounds in \eqref{eq:growth2} are obvious and that the second upper bound follows from the first one and \eqref{DB}. Hence, it suffices to establish the first upper bound.
	
	Let $f$ be the generating function for the offspring probabilities given by \eqref{branchingprob},
	\bb\label{f1}
	f(x) = \sum_{n=0}^\infty p(n)x^n = \frac{1-X(g_c(k))}{1-xX(g_c(k))} = \frac{1}{1-m(x-1)}\,,
	\ee
	where the last equality follows from \eqref{mean}, 
	and let $f_r$ denote the probability generating function for the size of the $r$th generation of the corresponding BGW process,
	$$
	f_r(x) = {\mathbb E}_\lambda (x^{|D_r|})\,,\quad r\geq 1\,,
	$$
	where $\lambda$ is the associated probability measure on $\cT$. It is a standard result (see e.g. \cite{athreyaney1972branching}) that
	\bb\label{itfr}
	f_{r+1} = f \circ f_r\quad\mbox{for $r\geq 2$} \quad\mbox{and }\;  f_2 =f\,.
	\ee
	Similarly, let $g_r$ be the probability generating function for the size of the $r$'th generation w.r.t.  the local limit measure $\nu^{(\mu)}$. Using the same notation as in the proof of Corollary \ref{growth1} for the corresponding single spine tree $\hat T$ and its branches grafted to the left and right along the spine, it follows from \eqref{decompDr1} that 
	\bb\label{gr}
	g_r(x) = {\mathbb E}_\mu(x^{|D_r|}) = x\prod_{k=2}^r f_k(x)^2\,.
	\ee
	
	We claim there exist constants $b>1$ and $0<c<1$ such that
	\bb\label{frbound}
	f_r(x) \leq 1 + c^{r-1}(x-1)\,,\quad 1\leq x\leq b\,,\;r\geq 2\,.
	\ee
	Indeed, by \eqref{f1} we have
	$$
	f(x) = 1+m(x-1) + \frac{m^2(x-1)^2}{1-m(x-1)}\,,
	$$
	and since $0<m<1$ we can choose $b>1$ such that
	$$
	0< \frac{m^2(b-1)}{1-m(b-1)} < 1-m\,,
	$$
	and hence 
	$$
	f(x) \leq 1 + c(x-1)\quad\mbox{for $1\leq x\leq b$,}
	$$
	where 
	$$
	c := m+ \frac{m^2(b-1)}{1-m(b-1)} <1\,.
	$$
	A simple inductive argument using \eqref{itfr} then implies \eqref{frbound}. 
	
	Combining \eqref{gr} and \eqref{frbound}, we obtain
	$$
	g_r(x) \leq x\prod_{k=2}^r e^{2c^{k-1}(x-1)} \leq b\cdot e^{2\frac{x-1}{1-c}}
	$$
	for $1\leq x\leq b$. By the Chebychev inequality 
	$$
	{\mathbb E}_\mu( e^{\theta |D_r|})\; \geq\; e^{\theta\lambda} \nu^{(\mu)}(\{|D_r|\geq \lambda\})\,,
	$$
	it therefore follows by choosing $\theta =\ln b$ that
	$$
	\nu^{(\mu)}(\{|D_r|\geq \lambda\}) \;\leq\; g_r(b)\cdot b^{-\lambda}\; \leq\; e^{2\frac{b-1}{1-c}} \cdot b^{1-\lambda}\,,
	$$
	for any $\lambda > 0$. Thus, choosing $\lambda = C_1\cdot\ln r$, we get that 
	$$
	\sum_{r=1}^\infty \nu^{(\mu)}(\{|D_r|\geq  C_1\cdot\ln r\})\; <\;\infty\,,
	$$
	provided $C_1>(\ln b)^{-1}$. Invoking the Borel-Cantelli lemma, the first upper bound in \eqref{eq:growth2} then holds for $r$ sufficiently large for a.e. $\hat T$, which concludes the proof.

	\subsection*{Proof of Theorem \ref{asvolgrowth}}
	
	\noindent\emph{Upper bound.} \; We  employ a method similar to the one of the previous proof. Note first that the generating function $f$ for the offspring probabilities \eqref{offspring2} (corresponding to $m=1$ in the previous proof) fulfills
	$$
	f(x) = \frac{1}{2-x} = x + (x-1)^2 + \frac{(x-1)^2}{2-x} \leq x + 3(x-1)^3\,,\quad\mbox{for $1\leq x\leq \frac 32$}.
	$$
	The probability generating function for $|D_r|$ w.r.t. the corresponding BGW process is then given by \eqref{itfr}. It is easily shown by induction that
	\bb\label{appest1}
	f_r(x) \leq x+12r(x-1)^2\,,\quad 0\leq x-1\leq \frac{1}{12(r-1)}\,.
	\ee
	Consider now $T\in\chi^{-1}(\cB_{\frac 1r}(T^s))$ for given $r\geq 1$ and $T^s\in \cT^s$. As previously, let $T_{(i,n)}$ denote the finite branch of $T$ grafted in sector $(i,n)$ of $T_0^s = B_r(T^s)$ and let $D^f_r(T)$ denote the set of vertices of finite type in $T$ at height $r$, i.e.
	$$
	|D^f_r(T)| = \sum_{(i,n)} |D_{r-h_i+1}(T_{(i,n)})|\,.
	$$
	Using that the branches $T_{(i,n)}$ are i.i.d. according to $\rho$ by Theorem \ref{backbonedist}, it follows that 
	$$
	\mathbb E_{\mu}(x^{|D_r^f|}) = \tilde{\mathbb E}_{\mu}\Big(\prod_{h=1}^{r-1} f_{r-h+1}(x)^{|D_h(T^s)|+|D_{h+1}(T^s)|}\Big)\,,
	$$
	since the number of sectors associated to $T^s$ at height $h$ equals $|D_h(T^s)|+|D_{h+1}(T^s)|$. 
	
	Writing $|D_h|= 1+\zeta_1+\zeta_2+\dots+\zeta_{h-1}$ on $\cT^s$,  then gives
	\begin{align}
		\mathbb E_{\mu}(x^{|D_r^f|})& = \Big(\prod_{k=2}^r f_k(x)^2\Big) \tilde{\mathbb E}_{\mu}\Big(\prod_{l=1}^{r-2}\Big( f(x) \prod_{m =l+1}^{r-1}f_{r-m+2}(x)f_{r-m+1}(x)\Big)^{\zeta_l}f(x)^{\zeta_{r-1}}\Big)\nonumber\\ &= \Big(\prod_{k=2}^r f_k(x)^2\Big) \tilde{\mathbb E}_{\mu}\Big(\prod_{l=1}^{r-2}\Big(f_{r-l+1}(x) \prod_{m =2}^{r-l}f_{m}(x)^2\Big)^{\zeta_l}\cdot f(x)^{\zeta_{r-1}}\Big)\,.\nonumber 
	\end{align}
	Here, the right-hand side can be computed by use of  Corollary \ref{cor:nutilde} and
	$$
	\tilde{ \mathbb E}_{\mu}(t^{\zeta_k}) = e^{\mu(t-1)}\,,\quad t\in\mathbb R\,,
	$$
	yielding the result 
	\bb\label{nutildegenfct}
	\mathbb E_{\mu}(x^{|D_r^f|}) = e^{\mu(f(x) -1)} \Big(\prod_{k=2}^r f_k(x)^2\Big) \prod_{l=1}^{r-2}\,e^{\mu \Big(f_{r-l+1}(x) \prod_{m =2}^{r-l}f_{m}(x)^2-1\Big)}\,.
	\ee
	Setting $x=e^\theta$ and using Chebychev's inequality, we have
	\bb\label{appcheb}
	\nu^{(\mu)}(\{|D^f_r|\geq \lambda r^2\}) \leq e^{-\theta\lambda r^2} \mathbb E_{\mu}(e^{\theta |D_r^f|}) \,,
	\ee
	for $\theta, \lambda>0$. From this, a useful bound is obtained by choosing  $\theta = r^{-2}$ and bounding the factors on the right-hand side of  \eqref{nutildegenfct} as follows. Noting that \eqref{appest1} implies 
	\bb\label{appest2}
	f_k(x) \leq  1+ 2(x-1) \leq e^{2(x-1)}\,,\quad \mbox{for $0\leq x-1\leq \frac{1}{12r}$ and $2\leq k\leq r$} \,,
	\ee
	we  get
	$$
	f_k(e^{r^{-2}}) \leq e^{4 r^{-2}}\quad\mbox{for $2\leq k\leq r$ and $r\geq 24$.} 
	$$
	Hence, 
	$$
	f_{r-l+1}(e^{r^{-2}}) \prod_{m =2}^{r-l}f_{m}(e^{r^{-2}})^2 -1 \leq e^{8 r^{-1}} -1 \leq 16 r^{-1} \,,\quad 1\leq l\leq r-2\,, 
	$$
	and consequently we obtain for the last product in \eqref{nutildegenfct} the bound
	$$
	\prod_{l=1}^{r-2}\,e^{\mu \Big(f_{r-l+1}(e^{r^{-2}}) \prod_{m =2}^{r-l}f_{m}(e^{r^{-2}})^2-1\Big)} \leq e^{16\mu}\,,
	$$
	for $r\geq 24$. Since these estimates likewise show that the first two factors in \eqref{nutildegenfct} are  bounded by constants for $x=e^{r^{-2}}$, it follows from  \eqref{nutildegenfct} and \eqref{appcheb} that 
	$$ 
	\nu^{(\mu)}(\{|D^f_r|\geq \lambda r^2\})\leq  C\cdot e^{-\lambda}\,,  
	$$
	where $C>0$ is a constant independent of $r$. 
	
	Choosing $\lambda = a\log r$, where $a>1$, this implies that
	$$
	\sum_{r=1}^\infty \nu^{(\mu)}(\{|D^f_r|\geq a r^2\log r\}) < \infty \,.
	$$
	Hence, by the Borel-Cantelli lemma, there exists $r_0(T)\in\mathbb N$ for $\nu^{(\mu)}$-a.e. $T$ such that 
	$$
	|D^f_r(T)| \leq ar^2\ln r\quad \mbox{for all $r\geq r_0(T)$}\,.
	$$
	Since $|D_r(T)| = |D_r^f(T)| + |D_r(\chi(T)|$, it follows from Corollary \ref{expectedbackbone} iii) and  \eqref{DB}  that the upper bound in \eqref{ashaus} holds for any $C''_2>\frac 13$.
	
	\bigskip
	
	\noindent\emph{Lower bound.}\;  For $N\in\mathbb N$, let 
	$$
	{\mathcal A}_N^s = \{T^s\in \cT^s\mid  \frac 12 \mu (r+1)\leq |D_r(T^s)|\;\leq\; 3\mu r\quad\mbox{for all $r\geq N$}\}.
	$$
	Then ${\mathcal A}_1^s\subseteq {\mathcal A}_2^s\subseteq {\mathcal A}_3^s\subseteq\dots$ and 
	\bb\label{ascend}
	\tilde\nu^{(\mu)}({\mathcal A}_N^s) \to 1\quad\mbox{as $N\to\infty$}
	\ee
	as a consequence of Corollary \ref{expectedbackbone} iii).  
	
	Let $N$ be fixed, as well as $r\geq 3N$. Given $T^s\in {\mathcal A}_N^s$, set $T_0^s = B_r(T^s)$ and note that 
	\bb\label{Tscond}
	\frac 12\mu\Big\lfloor\frac  r3\Big\rfloor\leq \big|D_{\lfloor\frac r3\rfloor}(T_0^s)\big|\leq \mu r\,.
	\ee
	Since $T_0^s$ has no leaves at heights less than $r$ there exist disjoint paths $\omega_i,\; i\in D_{\lfloor\frac r3\rfloor}(T_0^s)$, in $T_0^s$ of length $\lfloor \frac r3\rfloor$ such that $\omega_i$ connects $i$ to a vertex $j_i\in D_{2\lfloor\frac 3 r\rfloor}(T_0^s)$. Let us denote the vertices in $\omega_i$ by $v_0^i, v_1^i,\dots, v_{\lfloor \frac r3\rfloor}^i$, such that $v_0^i =i$ and $v_{\lfloor \frac r3\rfloor}^i  =j_i$. Now, consider $T\in\chi^{-1}(\cB_{\frac 1r}^s(T_0^s))$ and denote, for each $i\in D_{\lfloor\frac r3\rfloor}(T_0^s)$, by $T_k^i,\;0 \leq k\leq \lfloor \frac r3\rfloor$, the branches of $T$ grafted in sectors $(v_k^i,1)$ along $\omega_i$, respectively.  With respect to  $\nu^{(\mu)}$, conditioned on the set $ \chi^{-1}(\cB_{\frac 1r}^s(T_0^s))$, the branches $T_k^i$ are i.i.d. according to $\rho$, 
	by Theorem \eqref{backbonedist} ii). Hence, setting 
	\bb\label{defq}
	q(\lambda) = \nu^{(\mu)}(\{\sum_{k=0}^{\lfloor\frac r3\rfloor}|B_{\big\lfloor\frac r3\big\rfloor}(T_k^i)| \leq \lambda \big(\big\lfloor\frac r3\big\rfloor\big)^2\} \mid  \chi^{-1}(\cB_{\frac 1r}^s(T_0^s)))
	\ee
	for $\lambda >0$, we have that $q(\lambda)$ is independent of  $i\in D_{\lfloor\frac r3\rfloor}(T_0^s)$ and  that
	\bb\label{sumbound1}
	q(\lambda) \leq \Big(\rho(\{|B_{\lfloor\frac r3\rfloor}|\leq  \lambda \big(\big\lfloor\frac r3\big\rfloor\big)^2\})\Big)^{ \lfloor\frac r3\rfloor}\,.
	\ee
	We next establish an upper bound on the last expression by making use of the fact that there exists a constant $c'>0$ such that
	\bb\label{rhobound}
	\rho(\{|B_R|\geq  R^2\}) \geq \frac{c'}{R}\,,\quad R\geq 1\,.
	\ee
	This is a consequence of the following estimates:
	\begin{align}
		\rho(\{|B_R|\geq  R^2\}) &= 2\sum_{\substack{|T|\geq R^2\\ h(T)<R}} 4^{-|T|} + 2\sum_{\substack{|T|\geq  R^2\\ h(T)=R}} 4^{-|T|} 2^{|D_R(T)|}\nonumber\\ 
		&\geq 2\sum_{\substack{|T|\geq R^2\\ h(T)\leq R}} 4^{-|T|} \nonumber \\ 
		&= 2\sum_{N\geq  R^2} A_{R,N} 4^{-N} \nonumber\\
		& \geq 2\sum_{N\geq R^2} \frac{1}{R+1}\tan^2\frac{\pi}{R+1} \big(1+\tan^2\frac{\pi}{R+1}\big)^{-N}\nonumber\\
		& \geq \frac{2}{R+1} \sum_{N\geq R^2} \frac{\pi^2}{(R+1)^2} \big(1+\tan^2\frac{\pi}{R+1}\big)^{-N}\,, 
	\end{align}
	where we have used \eqref{rhoball} and \eqref{eq:A}, and we note that the last sum is clearly finite for all $R\geq 1$ and converges to $\int_{\pi^2}^\infty e^{-x}dx = e^{-\pi^2}$ as $R\to\infty$, thus proving \eqref{rhobound}. 
	
	For $ 0<\lambda <1$ it then follows that
	$$
	\rho(\{|B_R|\leq \lambda R^2\}) \leq \rho(\{|B_{\lfloor\sqrt\lambda R\rfloor+1}|\leq \lambda R^2\}) \leq 1-\frac{c'}{\sqrt\lambda R+1}
	$$
	for $R\geq 1$, and hence \eqref{sumbound1} implies 
	\bb\label{sumbound2}
	q(\lambda) \leq e^{-\frac{c'}{2\sqrt\lambda}}\,,
	\ee
	if $\lfloor\frac r3\rfloor\geq \lambda^{-\frac  12}$.
	
	Returning to $T\in\chi^{-1}(\cB_{\frac 1r}^s(T_0^s))$, we have 
	\begin{align}
		&~~~~~~ |B_ r(T)| \leq \frac{\lambda}{4} {\Big\lfloor\frac r3\Big\rfloor}^3 \nonumber\\
		&\Rightarrow~~ \sum_{i\in D_{\lfloor\frac r3\rfloor}(T_0^s)} \sum_{k=0}^{\lfloor \frac r3\rfloor} |B_{\lfloor\frac r3\rfloor}(T_k^i)|\; \leq \;\frac{\lambda}{4} {\Big\lfloor\frac r3\Big\rfloor}^3\nonumber \\
		& \Rightarrow~~ \Big|\Big\{i\in D_{\lfloor\frac r3\rfloor} (T_0^s)\,\Big|\; \sum_{k=0}^{\lfloor \frac r3\rfloor} |B_{\lfloor\frac r3\rfloor}(T_k^i)|\; \leq\; \frac{\lambda}{\mu} {\Big\lfloor\frac r3\Big\rfloor}^2\Big\}\Big| \;\geq\; \frac{\mu}{4}{\Big\lfloor\frac r3\Big\rfloor}\,,
	\end{align} 
	where the last implication follows from the first inequality of \eqref{Tscond}. Hence, it follows from the independence of the events $\{|B_{\lfloor\frac r3\rfloor}(T_k^i)| \leq \lambda \big(\lfloor\frac r3\rfloor\big)^2\},\, i\in D_{\lfloor\frac r3\rfloor}(T_0^s)|\}$, and the definition \eqref{defq} of $q$,  that 
	\begin{align}
		\nu^{(\mu)}\Big(\Big\{ |B_{ r}|& \leq \frac{\lambda}{4} {\Big\lfloor\frac r3\Big\rfloor}^3\Big\}\,\Big|\; \chi^{-1}(\cB_{\frac 1r}^s(T_0^s))\Big) \nonumber\\
		&\leq \sum_{m\geq \frac 14 \mu {\lfloor\frac r3\rfloor}}\binom{|D_{\lfloor\frac r3\rfloor}(T_0^s)|}{m} q(\lambda\mu^{-1})^m\big(1-q(\lambda\mu^{-1})\big)^{|D_{\lfloor\frac r3\rfloor}(T_0^s)|-m}\nonumber\\
		& \leq  \sum_{m\geq \frac 14 \mu {\lfloor\frac r3\rfloor}}\binom{\lfloor\mu r\rfloor}{m} q(\lambda\mu^{-1})^m \leq \mbox{cst}\cdot q(\lambda\mu^{-1})^{\frac 1{12} \mu r}\cdot 2^{\mu r}\,,
	\end{align}
	where the upper bound in \eqref{Tscond} has also been used. Recalling \eqref{sumbound2}, we may choose $\lambda=\lambda(\mu)$ small enough, independently of $r$ and $T_0^s$, such that $q_0 :=2q\Big(\frac{\lambda(\mu)}{\mu}\Big)^{\frac{1}{12}}<1$, yielding  
	\bb\label{sumbound3}
	\nu^{(\mu)}\Big(\Big\{ |B_{ r}| \leq \frac{\lambda(\mu)}{4} {\Big\lfloor\frac r3\Big\rfloor}^3\Big\}\,\Big|\; \chi^{-1}(\cB_{\frac 1r}^s(T_0^s))\Big)\; \leq\;  \mbox{cst}\cdot q_0^{\mu r}\,,
	\ee
	provided $\lfloor\frac r3\rfloor \geq \mbox{max}\{N, \Big(\frac{\mu}{\lambda(\mu)}\Big)^{\frac 12}\}$. Since this holds uniformly in $T^s\in{\mathcal A}_N^s$, with $T_0^s = B_r(T^s)$, for $r$ large enough, we conclude for such $r$ that  
	$$
	\nu^{(\mu)}\Big(\Big\{ |B_{ r}| \leq \frac{\lambda(\mu)}{4} {\Big\lfloor\frac r3\Big\rfloor}^3\Big\}\cap \chi^{-1}({\mathcal A}_N^s)\Big) \;\leq\; \mbox{cst}\cdot q_0^{\mu r}\,.
	$$
	In particular, it follows that 
	$$
	\sum_{r=1}^\infty \nu^{(\mu)}\Big(\Big\{ |B_{ r}| \leq \frac{\lambda(\mu)}{4} {\Big\lfloor\frac r3\Big\rfloor}^3\Big\}\cap \chi^{-1}({\mathcal A}_N^s)\Big) < \infty\,.
	$$
	Hence, the Borel-Cantelli lemma implies that $\{|B_{ r}| \leq \frac{\lambda(\mu)}{4} {\lfloor\frac r3\rfloor}^3\;\mbox{i.o.}\;\}\,\cap\, \chi^{-1}({\mathcal A}_N^s)$ is a nullset. Since this holds for arbitrary $N$, it follows from \eqref{ascend} that $\{ |B_{ r}| \leq \frac{\lambda(\mu)}{4} {\lfloor\frac r3\rfloor}^3\;\mbox{ i.o.} \,\}$ is a nullset, which hence completes the proof of the lower bound in \eqref{ashaus} with $0< C''_1 < \frac{\lambda(\mu)}{108}$.

	\bibliographystyle{abbrv}

\end{document}